\documentclass[journal,10 pt,oneside,web]{ieeecolor} %

\usepackage{graphicx}
\usepackage{textcomp}
\usepackage{epsfig} 
\usepackage{epsfig,color,amsmath,cite}

\usepackage{wrapfig}
\usepackage{bm}
\usepackage{epstopdf}
\usepackage{amssymb}
\usepackage{url}

\usepackage{multirow}
\usepackage{hhline}
\usepackage{booktabs}
\usepackage[linesnumbered,boxed,commentsnumbered,ruled,vlined,longend]{algorithm2e}
\usepackage{comment}

\DeclareMathOperator{\trace}{tr}
\DeclareMathOperator{\rank}{rank}

\DeclareMathOperator*{\argmax}{arg\ max}

\makeatother
\DeclareMathAlphabet\mathbfcal{OMS}{cmsy}{b}{n}



\usepackage{stackengine}

\newcommand{\prob}{\mathbb{P}}

\newcommand{\m}{\boldsymbol}
\allowdisplaybreaks[4]
\makeatletter
\let\NAT@parse\undefined
\makeatother

\usepackage[colorlinks = true,
linkcolor = blue,
urlcolor  = black,
citecolor = blue,
anchorcolor = blue]{hyperref}

\newcommand{\mc}[1]{\mathcal{#1}}
\newcommand{\mbb}[1]{\mathbb{#1}}
\newcommand{\mr}[1]{\mathrm{#1}}
\usepackage[framemethod=TikZ]{mdframed}
\mdfdefinestyle{MyFrame}{%
	linecolor=black,
	outerlinewidth=1.25pt,
	roundcorner=1.25pt,
	innerrightmargin=5pt,
	innerleftmargin=5pt,}
	
\usepackage[noabbrev]{cleveref}

\usepackage{mathtools}

\DeclarePairedDelimiter\abs{\lvert}{\rvert}%
\DeclarePairedDelimiter\norm{\lVert}{\rVert}%

\makeatletter
\let\oldabs\abs
\def\abs{\@ifstar{\oldabs}{\oldabs*}}
\let\oldnorm\norm
\def\norm{\@ifstar{\oldnorm}{\oldnorm*}}
\makeatother

\usepackage[english]{babel}
\usepackage[utf8]{inputenc}
\usepackage[super]{nth}

\usepackage{graphicx}
\usepackage{float}
\usepackage[caption = false]{subfig}

\usepackage{array}
\usepackage{threeparttable}

\usepackage[english]{babel}
\usepackage[utf8]{inputenc}
\usepackage[super]{nth}

\usepackage{pifont}

\usepackage{mathtools}

\usepackage{mleftright}
\usepackage{xparse}

\NewDocumentCommand{\evaluat}{sO{\big}mm}{%
	\IfBooleanTF{#1}
	{\mleft. #3 \mright|_{#4}}
	{#3#2|_{#4}}%
}

\hypersetup{
pdfauthor={Mohamad H. Kazma},
pdftitle={Connections between Gramians and DPPs}
}
\let\labelindent\relax
\usepackage{enumitem}
\setlist[itemize]{leftmargin=*}

\usepackage{lettrine} 

\usepackage{balance} 

\usepackage{colortbl}

\usepackage{pgfplots}
\pgfplotsset{compat=1.18}

\usepackage{tikz}
\usetikzlibrary{shapes.geometric, arrows.meta, decorations.pathmorphing, decorations.markings, calc, 3d}
\usetikzlibrary{positioning, fit, backgrounds}
\usetikzlibrary{arrows.meta}  
\usetikzlibrary{positioning}
\usetikzlibrary{fit}		
\usetikzlibrary{decorations.pathreplacing, calc}
\usetikzlibrary{tikzmark, calc} 
\tikzset{
	tangent/.style={decoration={
			markings,mark=at position #1 with {
				\coordinate (ta) at (0,0);
				\coordinate (tb) at (0.1,0);
			}
		},postaction=decorate},
	tangent/.default=0.5
}

\definecolor{Vgray}{RGB}{247,247,247}
\definecolor{Vedge}{RGB}{200,200,200}

\definecolor{Sub1purple}{RGB}{200,200,255}
\definecolor{Sub1edge}{RGB}{205,215,238}

\definecolor{Sub2yellow}{RGB}{251,234,206}
\definecolor{Sub2edge}{RGB}{251,234,206} 

\definecolor{Sub3green}{RGB}{215,238,210}
\definecolor{Sub3edge}{RGB}{215,232,210}

\definecolor{Sub1color}{rgb}{0.85, 0.9, 1}
\definecolor{Sub2color}{rgb}{1, 0.9, 0.7}
\definecolor{Sub3color}{rgb}{0.8, 1, 0.8}
\definecolor{Sub4color}{RGB}{200,200,255}

\definecolor{Sub4edge}{RGB}{205,215,238}

\definecolor{Sub5color}{RGB}{255,210,210}
\definecolor{Sub5edge}{RGB}{240,180,180}

\definecolor{Sub5edge}{RGB}{140,120,80}
\definecolor{Sub5color}{RGB}{255,245,225}

\definecolor{Sub6edge}{RGB}{100,130,150}
\definecolor{Sub6color}{RGB}{220,240,255}

\definecolor{Sub7edge}{RGB}{130,100,150}
\definecolor{Sub7color}{RGB}{245,230,255}

\definecolor{Sub8edge}{RGB}{150,130,100}
\definecolor{Sub8color}{RGB}{255,250,230}

\definecolor{Sub9edge}{RGB}{90,140,140}
\definecolor{Sub9color}{RGB}{220,250,250}

\definecolor{Sub10edge}{RGB}{110,110,110}
\definecolor{Sub10color}{RGB}{240,240,240}

\tikzset{
	highlightP1/.style={
		fill=Sub2yellow, rounded corners=2pt, inner sep=2pt, anchor=base},
	highlightP2/.style={
		fill=Sub3green, rounded corners=2pt, inner sep=2pt, anchor=base}
}

\definecolor{LightBlue}{rgb}{0.7, 0.85, 1}
\definecolor{LightOrange}{rgb}{1, 0.8, 0.6}
\definecolor{LightGreen}{rgb}{0.843, 0.933, 0.824}

\tikzstyle{box} = [rectangle, draw, thick, rounded corners, 
text centered, font=\footnotesize, minimum height=1cm, text width=4.5cm]
\tikzstyle{colorbox1} = [rectangle, draw=black, rounded corners, fill=Sub1purple, 
text centered, font=\footnotesize, minimum height=1cm, text width=4.5cm]
\tikzstyle{colorbox2} = [rectangle, draw=black, rounded corners, fill=Sub2yellow, 
text centered, font=\footnotesize, minimum height=1cm, text width=4.5cm]
\tikzstyle{colorbox3} = [rectangle, draw=black, rounded corners, fill=Sub3green, 
text centered, font=\footnotesize, minimum height=1cm, text width=4.5cm]
\tikzstyle{colorbox4} = [rectangle, draw=black, rounded corners, fill={rgb,255:red,247; green,247; blue,247},text centered, font=\scriptsize, minimum height=0.5cm, text width=3cm]

\tikzstyle{arrow} = [thick,->,>=stealth]
\tikzstyle{dashedarrow} = [thick, dashed, ->, >=stealth]
\usepackage{listings}
\usepackage{circuitikz}

\usepackage{makecell}

\usepackage{flowchart}

\newcommand{\introstart}[2]{\noindent \lettrine[lines=2]{#1}{#2}}

\definecolor{LCSS}{HTML}{004498}

\definecolor{FigBlue}{HTML}{4477BB}
\definecolor{FigGreen}{HTML}{66AA66}
\definecolor{FigPurple}{HTML}{AA88CC}
\definecolor{FigGold}{HTML}{CCAA44}

\newcommand{\parstartc}[1]{\noindent \textbf{\textcolor{LCSS}{#1.}}\;}

\newcommand{\Rn}[1]{\mathbb{R}^{#1}}

\newcommand{\logdet}[1]{\mr{log} \mr{det}(#1)}

\newcommand{\titlesc}[1]{\title{\Large \vspace{0.7cm} \LARGE \centering {\textsc{{#1}}}}}

\renewenvironment{proof}{%
	\par\vspace{0.3em}%
	\noindent\textbf{\textit{\textcolor{LCSS}{Proof.}}}\enspace%
}{%
	\hfill$\blacksquare$%
	\par\vspace{0.3em}%
}
\newtheorem{theorem}{\textbf{\textcolor{LCSS}{Theorem}}}
\newtheorem{mylem}{\textbf{\textcolor{LCSS}{Lemma}}}
\newtheorem{mydef}{\textbf{\textcolor{LCSS}{Definition}}}

\newtheorem{mycor}{\textbf{\textcolor{LCSS}{Corollary}}}
\newtheorem{myprs}{\textbf{\textcolor{LCSS}{Proposition}}}
\newtheorem{exmpl}{\textbf{\textcolor{LCSS}{Example}}}
\newtheorem{problem}{\textbf{\textcolor{LCSS}{Problem}}}

\makeatletter
\def\thm@space@setup{%
  \thm@preskip=8pt plus 1pt minus 1pt%
  \thm@postskip=8pt plus 1pt minus 1pt%
}
\def\@begintheorem#1#2{%
	\item[\hskip\labelsep\textbf{\textcolor{LCSS}{\textit{#1 #2}}}]%
	\@ifnextchar[{\@withinfo}{\textbf{\textcolor{LCSS}{.}}\enspace\ignorespaces}}
	
	\def\@withinfo[#1]{
\textbf{\textcolor{LCSS}{\textit{ (#1)}}}\textcolor{LCSS}{.} \ignorespaces}
\makeatother

\makeatletter
\def\@thm#1#2{%
  \refstepcounter{#1}%
  \@ifnextchar[{\@ythm{#1}{#2}}{\@xthm{#1}{#2}}}
\def\@xthm#1#2{%
  \@begintheorem{#2}{\csname the#1\endcsname}\ignorespaces}
\def\@ythm#1#2[#3]{%
  \@opargbegintheorem{#2}{\csname the#1\endcsname}{#3}\ignorespaces}
\def\@endtheorem{\par}
\def\@begintheorem#1#2{%
	\par\noindent\textbf{\textit{\textcolor{LCSS}{#1 #2.}}}%
	\@ifnextchar[{\@withinfo}{\enspace\ignorespaces}}
\def\@withinfo[#1]{
	\textbf{\textcolor{LCSS}{ (#1)}}\enspace\ignorespaces}
\def\@opargbegintheorem#1#2#3{%
	\par\noindent\textbf{\textit{\textcolor{LCSS}{#1 #2. (#3).}}}%
	\enspace\ignorespaces}
\makeatother

\newcommand{\logdett}{\mr{log} \mr{det}}

\hypersetup{linkcolor=LCSS,citecolor=LCSS}

\newcommand{\theoref}[1]{\hyperref[#1]{Theorem~\ref*{#1}}}
\newcommand{\lemref}[1]{\hyperref[#1]{Lemma~\ref*{#1}}}
\newcommand{\coref}[1]{\hyperref[#1]{Corollary~\ref*{#1}}}
\newcommand{\propref}[1]{\hyperref[#1]{Proposition~\ref*{#1}}}
\newcommand{\defref}[1]{\hyperref[#1]{Definition~\ref*{#1}}}
\newcommand{\remref}[1]{\hyperref[#1]{Remark~\ref*{#1}}}
\newcommand{\exmpref}[1]{\hyperref[#1]{Example~\ref*{#1}}}
\newcommand{\asmpref}[1]{\hyperref[#1]{Assumption~\ref*{#1}}}
\newcommand{\secref}[1]{\hyperref[#1]{Section~\ref*{#1}}}
\newcommand{\subsecref}[1]{\hyperref[#1]{Section~\ref*{#1}}}
\renewcommand{\eqref}[1]{\hyperref[#1]{(\ref*{#1})}}
\newcommand{\figref}[1]{\hyperref[#1]{Fig.~\ref*{#1}}}
\newcommand{\tabref}[1]{\hyperref[#1]{Table~\ref*{#1}}}
\newcommand{\probref}[1]{\hyperref[#1]{Problem~\ref*{#1}}}
\newcommand{\algref}[1]{\hyperref[#1]{Algorithm~\ref*{#1}}}
\newcommand{\appref}[1]{\hyperref[#1]{Appendix~\ref*{#1}}}

\usepackage{lcsys} %

\colorlet{subsectioncolor}{.} %
\definecolor{subsectioncolor}{RGB}{39,94,77} %

\titlesc{Connections Between Determinantal Point Processes and Gramians in Control}

\author{Mohamad H. Kazma$^{\diamond}$, \textit{Member, IEEE} and Ahmad F. Taha, \textit{Member, IEEE}  \vspace{-0.1cm}
	\thanks{
		$^\diamond$Corresponding author. This work is supported by National Science Foundation under Grant 2152450. The authors are with the Civil $\&$ Environmental Engineering and Electrical $\&$ Computer Engineering Departments, Vanderbilt University, 2201 West End Ave, Nashville, Tennessee 37235. Emails: mohamad.h.kazma@vanderbilt.edu, ahmad.taha@vanderbilt.edu.} 
}

\begin{document}
\newdimen\origiwspc%
\newdimen\origiwstr%
\origiwspc=\fontdimen2\font%
\origiwstr=\fontdimen3\font%

\maketitle
\thispagestyle{plain}
\pagestyle{plain}
\begin{abstract}
Determinantal point processes (DPPs) are probability models over subsets of a ground set that favor diverse selections while suppressing redundancy. That is, they tend to assign higher likelihood to collections where the elements complement one another instead of repeating the same information. This paper establishes connections between DPPs and control theory. By showing that the observability (controllability) Gramian of a linear dynamic system is a DPP over the sensor (control) node set, we provide a probabilistic and spectral perspective on sensor (actuator) selection for linear dynamic systems. This notion of probability here does \textit{not} represent stochastic uncertainty in the system dynamics; it instead represents a likelihood measure over sensor (actuator) configurations induced by the Gramian. To that end, we derive an effective observable rank condition, characterize the balance between individual node contributions and diversity, and establish node inclusion monotonicity and negative dependence properties. Finally, we show that the formulation recovers greedy optimization guarantees and admits a maximum a posteriori interpretation of the node selection problem. Numerical case studies on several network topologies corroborate the theoretical results.
\end{abstract}
\vspace{-0.5cm}
\section{Introduction and Contributions}\label{sec:Intro}
\introstart{O}{bservability} and controllability are fundamental system-theoretic properties~\cite{Kalman1963} that determine whether the internal state of a dynamic system can be reconstructed from output measurements or steered through actuators. For large scale systems, deploying sensors or actuators at every node is often infeasible due to sensing cost, communication constraints, and physical deployment limitations. The \emph{node selection problem}, equivalently the sensor or actuator selection problem, involves selecting a subset of nodes from a candidate set to optimize an observability or controllability Gramian measure; it therefore governs whether and to what extent a state variable can be reconstructed or actuated. The literature on this topic is inveterate and has been studied from several perspectives; we briefly discuss some of the literature next.

While the aforementioned problem is NP-hard, tractable approaches provide approximation guarantees. Convex relaxation approaches have been developed for sensor selection under D-optimality~\cite{Joshi2009}, and a convex formulation for placement under correlated observations is presented in~\cite{Patan2025}. Metrics characterizing tradeoffs between control energy and the number of control nodes are introduced in~\cite{Pasqualetti2014a}. Approximations that exploit the submodularity properties of Gramian measures are introduced in~\cite{Shamaiah2010,Summers2016,Tzoumas2016}, which yield greedy algorithms with a $(1 - 1/e)$-approximation guarantee~\cite{Nemhauser1978}. Greedy algorithms with matroid constraints have been applied to actuator placement~\cite{Guo2021}. Supermodular and weak submodular formulations have also been studied for trace and minimum eigenvalue Gramian measures~\cite{Clark2014,Sato2024}. Sensors and actuators are selected from a reduced order representation of the system through balanced model reduction in~\cite{Manohar2022}. Similar formulations address control node selection, where target controllability of complex networks is studied through greedy optimization~\cite{Ding2025}, and the controllability Gramian spectrum is computed and optimized for networked Laplacian systems with limited control placement in~\cite{Cao2024}.

\begin{figure}[t]
\vspace{-0.3cm}	
\centering
\includegraphics[width=0.95\columnwidth]{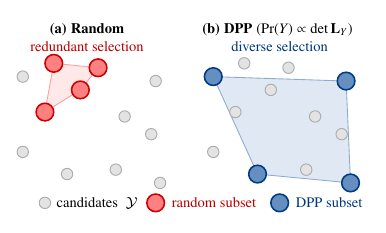}
\vspace{-0.3cm}
\caption{DPP versus independent (random) selection on $|\mc{Y}|=12$ candidates. Random sampling may produce redundant subsets while a DPP favors subsets with near orthogonal columns in $\m{L}$, thus suppressing redundant selections.} 
\label{fig:dpp_illustration}
\vspace{-0.15cm}
\end{figure}

Furthermore, randomized algorithms can reduce the computational burden of the aforementioned greedy algorithms while retaining performance guarantees. This includes actuator scheduling under Gramian metrics that are monotone and homogeneous~\cite{Siami2021}, greedy selection under weak submodularity~\cite{Hashemi2021,Kaya2025}, and uniform sensor sampling with probabilistic guarantees on spectral metrics of the nonlinear observability Gramian~\cite{Bopardikar2019}. Sensor selection has also been studied through the volume of the information ellipsoid of a configuration~\cite{Liu2023}, and through sparse sampling of network measurements~\cite{Amini2022}. Robust selections that retain utility across multiple submodular objectives are developed in~\cite{Kaya2024}. These approaches consider the selected configuration as a random subset drawn to approximate the selection objective; this is the perspective we adopt. In doing so, we exploit the determinantal structure of the Gramian and develop its probabilistic interpretation.

In this paper, we establish a connection between Gramians and a point process\footnote{A point process on a finite ground set $\mc{Y}$ is a probability measure on the power set $2^{\mc{Y}}$ that assigns a probability to each subset $Y \subseteq \mc{Y}$~\cite{Kulesza2012}.} that yields a probabilistic formulation of the node selection problem; it reveals spectral properties of the Gramian that are not apparent from existing formulations. A~\textit{determinantal point process} (DPP) is a probability distribution over subsets of a finite ground set in which each subset $Y \subseteq \mc{Y}$ is assigned a probability proportional to $\det(\m{L}_{Y})$, where $\m{L} \succeq \m{0}$ is a positive semidefinite matrix and $\m{L}_{Y}$ is the principal submatrix indexed by $Y$~\cite{Kulesza2012}. This determinantal structure assigns high probability to subsets for which the corresponding columns in $\m{L}$ are near orthogonal and suppresses subsets with redundant entries; see Fig.~\ref{fig:dpp_illustration}.

DPPs were introduced by Macchi~\cite{Macchi1975} to model fermion distributions at thermal equilibrium. DPPs were subsequently developed in random matrix theory~\cite{Hough2006,Tao2009}, probability~\cite{Lyons2003,Borodin2005}, and combinatorics~\cite{Kulesza2012}.
The determinantal structure yields tractable expressions for normalization, conditioning, and sampling~\cite{Kulesza2012}. In machine learning, DPPs have been applied to recommendation diversification and document summarization~\cite{Gillenwater2018,Kulesza2011}. For example, YouTube considers DPP ranking of search results to avoid showing users many near duplicate videos while maintaining relevance~\cite{wilhelm2018youtube}. Another example is in video summarization, where the objective is to select a small representative set of frames that cover the main content without repetition~\cite{zhang2016summary}. Furthermore, the column subset selection problem, where a matrix is approximated through a subset of its columns~\cite{Civril2012}, is addressed by sampling subsets with probability proportional to the squared determinant of the corresponding submatrix~\cite{Belhadji2020}. More recently, DPPs have been used for subset selection in wireless networks, where the kernel encodes both link quality and mutual information~\cite{Tu2025}. DPPs have also been used to model directionality in Gaussian data, where the kernel is derived from the covariance structure of the data rather than learned from samples~\cite{Ghosh2020a}. Greedy approximations of the maximum a posteriori (MAP) subset of a DPP under a cardinality constraint are studied in~\cite{Grosse2024}. %

While their applications have been extensive in computer science, DPPs have, to the best of our knowledge, \textit{not} been considered in the context of control, and their relationship to Gramians in control has not been established. To that end, we show that when each state variable is a candidate node within a network, the observability Gramian satisfies the conditions of a DPP and provides a probabilistic and spectral perspective on node selection. 
Note that the notion of probability herein is not that of stochastic uncertainty in the system; it is rather a measure of likelihood over node configurations and measurement subsets induced by the Gramian.

\parstartc{Contributions} The paper contributions are as follows.
\begin{itemize}[leftmargin=*]
	\item We show that when each state variable is a candidate node, the observability Gramian of a linear time-invariant (LTI) system satisfies the conditions of a DPP, such that each node configuration is assigned a probability proportional to the determinant of the corresponding Gramian principal submatrix. For a linear output matrix, we show that the determinant of the Gramian associated with a sensor subset is proportional to the probability that a DPP selects only rows of the observability matrix associated with those sensors. These results provide probabilistic interpretations of the node selection problem. For the latter, we show that the DPP probability of a sensor configuration equals the ratio between the determinant of the estimation error covariance of the full sensor set and that of the configuration.
	\item We show that the objective of node selection problem is proportional to the determinant of the error covariance of the remaining state variables, which relates the determinantal measure to the accuracy of state estimation. We derive several properties from the aforementioned results. Specifically, \textit{(i)} we introduce an effective observable rank that quantifies how many observability modes contribute for a given system; \textit{(ii)} we prove that improving system observability monotonically increases the inclusion probability of every node configuration; \textit{(iii)} we show that the complement of the selected node set characterizes the unobserved portion of the state-space, and \textit{(iv)} that the determinantal structure suppresses redundant node selection. These properties offer probabilistic insights for the node selection problem that are not captured by existing methods.
	\item We recover the submodularity of proposed measure using the Schur complement and reformulate it as a DPP MAP problem. This allows for greedy approximation guarantees while providing a distribution over configurations instead of one deterministic set. We validate the proposed framework on networks of varying topologies.
\end{itemize}

\parstartc{Broader impacts} The proposed framework returns a probabilistic characterization of node configurations rather than a single deterministic set. This characterization is relevant when the loss or corruption of measurements degrades observability~\cite{Babazadeh2022,Zhang2023c} and when sensing and communication budgets require scheduling the active nodes over the operation of the system~\cite{Vafaee2023,Abbas2018,Li2024f}. The ensemble of configurations obtained herein thereby provides alternatives that remain available when selected nodes fail. Such a perspective is essential for operating infrastructure networks under uncertainty, as it provides confidence intervals for estimation quality, risk-aware monitoring priorities, and reliability measures under changing operating conditions. This perspective also extends to actuator selection~\cite{Yan2015a,Abbas2024}, where a probabilistic characterization of the configurations would quantify how placement uncertainty affects the ability to regulate critical states, maintain safe operating limits, and recover from disturbances.

\begin{figure*}[t]
	\centering
	\scriptsize
	\resizebox{1\textwidth}{!}{
		\begin{tikzpicture}[
			node distance=0.45cm and 0.6cm,
			box/.style={draw=gray!70, rounded corners, fill=gray!8, text width=3.2cm, align=center},
			grouptitle/.style={font=\scriptsize\bfseries, align=center, text width=3.4cm},
			arrow/.style={->,>=stealth, thick, shorten <=4pt, shorten >=4pt}
			]

			\node[box, fill=gray!25] (thm1) {
				\textbf{Gramian as a DPP} \\[1pt]
				$\mc{P}_{\m{L}}(Y) \propto \det((\m{W}_o)_Y)$ \\
				\textit{\theoref{thm:gramian_dpp}}
			};

			\node[box, fill=gray!25, below=0.25cm of thm1] (thm2) {
				\textbf{Marginal kernel} \\[1pt]
				$\m{K} = \m{W}_o(\m{W}_o+\m{I})^{-1}$ \\ $n_{\mr{eff}} = \trace(\m{K})$ \\
				\textit{\theoref{thm:marginal_modal}}
			};

			\node[grouptitle, above=0.12cm of thm1] (t2) {Gramian DPP connection (\secref{sec:main_results})};

			\begin{scope}[on background layer]
				\node[draw=Vedge, fill=Vgray, rounded corners, inner sep=5pt, fit=(t2)(thm1)(thm2)] (P2) {};
			\end{scope}

			\node[box, left=of P2] (prob) {
				\textbf{Node selection problem (\secref{subsec:node_prob})} \\[1pt]
				$\max_{|Y| = r} \; \logdet{(\m{W}_o)_Y}$ \\
				\textit{\probref{prob:node_selection}}
			};

			\node[box, fill=gray!25, anchor=west] at ([xshift=0.75cm]P2.east |- thm1) (props) {
				\textbf{Inclusion properties} \\[1pt]
				$\prob(S \subseteq Y) = \det(\m{K}_S)$ \\
				\textit{\hyperref[cor:cond_incl]{Corollaries~\ref*{cor:cond_incl}-\ref*{cor:neg_dep}}}
			};

			\node[box, fill=gray!25, below=0.25cm of props] (map) {
				\textbf{Greedy MAP and $k$-DPP selection} \\[1pt]
				\textit{\lemref{lem:submod}, \propref{pro:kdpp_map}, \hyperref[alg:greedy_map]{Algorithms~\ref*{alg:greedy_map}-\ref*{alg:kdpp}}}
			};

			\node[grouptitle, above=0.12cm of props] (t3) {Resulting properties and methods (\secref{sec:node_selection})};

			\begin{scope}[on background layer]
				\node[draw=Vedge, fill=Vgray, rounded corners, inner sep=5pt, fit=(t3)(props)(map)] (P3) {};
			\end{scope}

			\node[box, right=of P3] (sim) {
				\textbf{Numerical case studies (\secref{sec:case_study})} \\[1pt]
				Probabilistic node selection, state-estimation, \\ node failures, actuator selection \\[1pt]
				\textit{$(\mr{Q}1)$-$(\mr{Q}4)$}
			};

			\draw[arrow] (prob) -- (P2);
			\draw[arrow] (P2.east) -- (P2.east -| P3.west);
			\draw[arrow] (P3) -- (sim);

		\end{tikzpicture}
	}
	\caption{Overview of the proposed framework. The node selection problem (\probref{prob:node_selection}) is recast through the observability Gramian as a DPP over the candidate sensor nodes when $\m{C} = \m{I}$ (\theoref{thm:gramian_dpp}), and the marginal kernel characterizes the spectral structure of the DPP and the effective observable rank (\theoref{thm:marginal_modal}). The connection yields the inclusion properties (\hyperref[cor:cond_incl]{Corollaries~\ref*{cor:cond_incl}--\ref*{cor:neg_dep}}) and the submodularity and MAP reformulation of \probref{prob:node_selection} that yield the greedy and $k$-DPP sampling algorithms (\lemref{lem:submod}, \propref{pro:kdpp_map}); these results are validated in~\secref{sec:case_study}.}
	\label{fig:framework}
	\vspace{0.2cm}
\end{figure*}

\parstartc{Notation} Let $\mbb{N}$, $\mbb{R}$, and $\mbb{R}_{>0}$ denote the sets of natural, real, and positive real numbers. Let $\Rn{n}$ and $\Rn{n \times m}$ denote the sets of real-valued column vectors of size $n$ and $n$-by-$m$ real matrices. The operators $\logdet{\m{A}}$, $\trace(\m{A})$, and $\rank(\m{A})$ return the log-determinant, trace, and rank of a square matrix $\m{A}$. For symmetric $\m{A}, \m{B}$, we write $\m{A} \succeq \m{0}$ for positive semidefiniteness, $\m{A} \succ \m{0}$ for positive definiteness, and $\m{A} \preceq \m{B}$ for $\m{B} - \m{A} \succeq \m{0}$. Calligraphic letters denote ground sets ($\mc{Y}$); italic capital letters ($Y$, $S$, $A$) denote subsets of ground sets. For a matrix $\m{M}$ and subset $S$, $\m{M}_S$ denotes the principal submatrix indexed by $S$. The symbol $\mc{P}_{\m{L}}$ denotes the probability measure on the power set $2^{\mc{Y}}$ induced by a DPP with L-ensemble kernel $\m{L}$, $\mc{P}_{\m{L}}(Y)$ denotes the probability assigned to a specific subset $Y$, and $Y \sim \mc{P}_{\m{L}}$ denotes a random subset drawn from this measure. We use $\prob(\cdot)$ for the probability of an event; $\prob(A \subseteq Y)$ denotes the probability that the subset $A$ is included in the random set $Y$. $\mbb{E}[|Y|]$ denotes the expected cardinality of the random subset $Y$, and $\mr{Var}[|Y|]$ denotes its variance.

\parstartc{Organization} The paper is organized as follows; see~\figref{fig:framework}. \secref{sec:prelim} presents the node selection problem and the DPP preliminaries. \secref{sec:main_results} establishes the Gramian DPP connection, the quality and diversity decomposition, and the marginal kernel with the effective observable rank. \secref{sec:node_selection} derives the inclusion properties and reformulates \probref{prob:node_selection} as a $k$-DPP MAP problem with greedy and sampling algorithms. Numerical results are presented in \secref{sec:case_study}, and \secref{sec:conclusion} concludes the paper. The DPP for linear output matrices is established in \appref{app:generalC}.

\section{Preliminaries and Problem Formulation}\label{sec:prelim}
\subsection{Sensor and actuator node selection problem}\label{subsec:lti}
Consider a discrete-time LTI system with state-space realization given as follows
\begin{equation}\label{eq:model_DT}
	\m{x}[{k+1}] = \m{A}\m{x}[k] + \m{B}\m{u}[k], \quad
	\m{y}[k] = \m{C}\m{x}[{k}],
\end{equation}
where $k \in \mbb{N}$ is the discrete-time index, $\m{x}[k] \in \Rn{n}$ is the state, $\m{u}[k] \in \Rn{m}$ is the input, and $\m{y}[k] \in \Rn{p}$ collects the outputs of $p$ candidate output nodes. The matrices $\m{A} \in \Rn{n \times n}$, $\m{B} \in \Rn{n \times m}$, and $\m{C} \in \Rn{p \times n}$ are the state-space matrices of the system. We denote the ground set of candidate output indices by $\mc{Y} = \{1, \ldots, p\}$ and the $i$-th row of $\m{C}$ by $\m{c}_i^\top \in \Rn{1 \times n}$, such that the $i$-th output node results in the scalar measurement equation $y_i[k] = \m{c}_i^\top \m{x}[k]$.

The system~\eqref{eq:model_DT} is said to be observable over a finite horizon $N \in \mbb{N}$ if and only if the observability matrix
\begin{equation}\label{eq:Obs_matrix}
	\m{O} := \left[\m{C}^\top,\; (\m{C}\m{A})^\top,\; \cdots,\; (\m{C}\m{A}^{N-1})^\top\right]^\top \in \Rn{Np \times n},
\end{equation}
satisfies $\rank(\m{O}) = n$. This observability rank condition is qualitative. A quantitative characterization is provided by the finite horizon observability Gramian
\begin{equation}\label{eq:LinObsGram}
	\m{W}_o := \sum_{k=0}^{N-1} (\m{A}^k)^\top \m{C}^\top \m{C}\, \m{A}^k = \m{O}^\top \m{O} \in \Rn{n \times n},
\end{equation}
which is positive semidefinite; it is positive definite if and only if the system is observable over horizon $N$~\cite{Sato2024}.

Based on the observability Gramian, common measures that quantify the observability of system~\eqref{eq:model_DT} include $\rank(\m{W}_o)$, $\lambda_{\min}(\m{W}_o)$, $\trace(\m{W}_o)$, and $\logdet{\m{W}_o}$~\cite{Pasqualetti2014a,Summers2016}. The $\trace(\m{W}_o)$ measure captures aggregate output sensitivity to state perturbations, while $\logdet{\m{W}_o}$ is related to D-optimal estimation through the volume of the state estimation error ellipsoid~\cite{Joshi2009}.  

For a subset $Y \subseteq \mc{Y}$ of output nodes, let $\m{C}_Y \in \Rn{|Y| \times n}$ be the submatrix of $\m{C}$ formed by the rows indexed by $Y$. The corresponding observability matrix is defined as
\begin{equation}\label{eq:obs_matrix_row}
	\m{O}_Y := \left[\m{C}_Y^\top,\; (\m{C}_Y\m{A})^\top,\; \cdots,\; (\m{C}_Y\m{A}^{N-1})^\top\right]^\top \in \Rn{N|Y| \times n}.
\end{equation}
The observability Gramian $\m{W}_o$, parameterized by the subset $Y$ of output nodes, can be written as 
\begin{equation}\label{eq:sub_gramian}
	\m{W}_o(Y) := \sum_{k=0}^{N-1} (\m{A}^k)^\top \m{C}_Y^\top \m{C}_Y\, \m{A}^k = \m{O}_Y^\top \m{O}_Y \in \Rn{n \times n}, 
\end{equation}
where $\m{W}_o(Y)$ quantifies the observability information contributed by the selected sensors. The sensor selection problem in~\cite{Joshi2009,Summers2016} maximizes Gramian measures of $\m{W}_o(Y)$.

\label{subsec:control}%
Furthermore, actuator node selection is characterized through the controllability Gramian. Each input node corresponds to a column of $\m{B}$ and for a subset $Y \subseteq \{1, \ldots, m\}$ of input nodes, let $\m{B}_Y \in \Rn{n \times |Y|}$ collect the columns of $\m{B}$ indexed by $Y$. For succinctness in notation, $Y$ denotes node sets for both sensor or input nodes. The controllability Gramian, parameterized by the subset $Y$, is written as
\begin{equation}\label{eq:LinConGram}
	\m{W}_c(Y) := \sum_{k=0}^{N-1} \m{A}^k \m{B}_Y \m{B}_Y^\top (\m{A}^k)^\top = \m{R}_Y \m{R}_Y^\top \in \Rn{n \times n},
\end{equation}
where $\m{R}_Y := [\m{B}_Y,\; \m{A}\m{B}_Y,\; \cdots,\; \m{A}^{N-1}\m{B}_Y] \in \Rn{n \times N|Y|}$ is the controllability matrix of the input nodes in $Y$. The system~\eqref{eq:model_DT} is controllable over the horizon $N$ if and only if $\m{W}_c := \m{W}_c(\{1, \ldots, m\})$ is positive definite. Note that the results established for sensor node selection extend to controllability under the mapping $(\m{A},\m{C}) \mapsto (\m{A}^\top,\m{B}^\top)$, such that $\m{W}_c = \m{W}_o(\m{A}^\top,\m{B}^\top)$. For brevity, we present the results from the perspective of observability. We validate the results by considering actuator selection in the case studies. 

\subsection{Node selection problem in discrete-time LTI systems}\label{subsec:node_prob}
We first consider the setting in which each state variable is a candidate sensor node, that is, $\m{C} = \m{I}$, where $p = n$ and $\mc{Y} = \{1, \ldots, n\}$. This is the standard setting in network controllability and observability studies~\cite{Liu2011,Pasqualetti2014a}. For a node subset $Y \subseteq \mc{Y}$, let $(\m{W}_o)_Y := [(\m{W}_o)_{ij}]_{i,j \in Y} \in \Rn{|Y| \times |Y|}$ denote the principal submatrix of the Gramian~\eqref{eq:LinObsGram} indexed by $Y$. Note that $(\m{W}_o)_{ij}$ is the inner product of the columns of $\m{O}$ associated with nodes $i$ and $j$, and $\det((\m{W}_o)_Y)$ is the squared volume of the parallelepiped spanned by the columns of $\m{O}$ indexed by $Y$. The node selection problem is formulated as follows.

\begin{problem}[Node selection problem]\label{prob:node_selection}
Given an LTI system~\eqref{eq:model_DT} with node set $\mc{Y}$ and a budget $r$, find the subset $Y \subseteq \mc{Y}$ of cardinality $r$ that maximizes the $\logdett$ of the principal submatrix of the Gramian~\eqref{eq:LinObsGram} as follows
\begin{equation}\label{eq:sns_problem}
	Y^\star = \argmax_{Y \subseteq \mc{Y},\, |Y| = r} \logdet{(\m{W}_{o})_Y}.
\end{equation}
\end{problem}

We note that the principal submatrix $(\m{W}_o)_Y$ is distinct from the Gramian $\m{W}_o(Y)$ in~\eqref{eq:sub_gramian}. Both are based on the observability matrix $\m{O}$, where the former restricts $\m{O}$ to the columns indexed by $Y$, while the latter enters through the output matrix $\m{C}_Y$, which restricts $\m{O}$ to the rows of the selected nodes in~\eqref{eq:obs_matrix_row}, and is the Gram matrix $\m{O}_Y^\top \m{O}_Y$ over all $n$ columns. The Cauchy-Binet formula~\cite[Sec.~0.8.7]{Horn1985} expands the determinant of this product as the sum of the squared $n \times n$ minors of $\m{O}_Y$ over its row subsets of size $n$. Thus, $\det(\m{W}_o(Y))$ is a sum of squared minors of $\m{O}_Y$ rather than a single principal minor of $\m{W}_o$. Furthermore, $\logdet{\m{W}_o(Y)} = -\infty$ whenever $\m{W}_o(Y)$ is rank deficient. In contrast, $\logdet{(\m{W}_o)_Y}$ is finite for every $Y$ whenever $\m{W}_o \succ \m{0}$, and thus \probref{prob:node_selection} remains well-posed. %

\probref{prob:node_selection} is combinatorial and NP-hard~\cite{Summers2016,Kulesza2012}. Greedy algorithms with $(1-1/e)$-approximation guarantees exist when the objective is monotone submodular\footnote{A set function $f: 2^{\mc{Y}} \to \mbb{R}$ is submodular if $f(A \cup \{j\}) - f(A) \geq f(B \cup \{j\}) - f(B)$ for all $A \subseteq B \subseteq \mc{Y}$ and $j \notin B$; it is monotone if $f(A) \leq f(B)$ whenever $A \subseteq B$~\cite{Nemhauser1978}.}~\cite{Nemhauser1978}, and submodularity of $\logdett$ Gramian objectives is established in~\cite{Shamaiah2010,Summers2016}. While the theoretical bound guarantees that the solution is at least $63\%$ of the optimal value, in practice this bound can be higher~\cite{Kazma2026a}. In the subsequent sections, we establish the submodularity of~\eqref{eq:sns_problem} in \lemref{lem:submod}, and the greedy guarantee for solving \probref{prob:node_selection} is thus recovered under the proposed approach. We note that when the total curvature\footnote{The total curvature of a monotone submodular function $f$ is defined as $c := 1 - \min_{j \in \mc{Y}} \frac{f(\mc{Y}) - f(\mc{Y} \setminus \{j\})}{f(\{j\})}$, where $c \in [0,1]$. $c = 0$ corresponds to a modular (linear) function and $c = 1$ to the worst case submodular function~\cite{Conforti1984,Vendrell2026}.} $c$ is known, the bound tightens to $\tfrac{1}{c}(1-e^{-c})$~\cite{Conforti1984}. Strategies that consider $c$ can further improve practical optimality gaps without increasing computational complexity~\cite{Vendrell2026}.

\subsection{Determinantal point processes}\label{subsec:dpp} 
To make this paper self-contained, we briefly introduce the definitions and properties of DPPs used in the subsequent sections. A DPP is a probability distribution over subsets of a ground set $\mc{Y} = \{1, \ldots, p\}$ in which the inclusion probabilities are governed by determinants of a positive semidefinite matrix~\cite{Kulesza2012}. The L-ensemble parameterization defines subset probabilities by principal minors of a positive semidefinite kernel \(\m{L}\), i.e., \(\mc{P}_{\m{L}}(Y)\propto \det(\m{L}_Y)\). The following definitions formalize the L-ensemble and its marginal kernel. 

\begin{mydef}[\textit{L-ensemble~\cite{Kulesza2012}}]\label{def:L_ensemble}
	Given a positive semidefinite matrix $\m{L} \in \Rn{p \times p}$, the \emph{L-ensemble} DPP defines a probability measure $\mc{P}_{\m{L}}$ for any subset $Y \subseteq \mc{Y}$ on $2^{\mc{Y}}$ by
	\begin{equation}\label{eq:dpp_L}
		\mc{P}_{\m{L}}(Y) :=\frac{\det(\m{L}_Y)}{\sum_{Y' \subseteq \mc{Y}} \det(\m{L}_{Y'})} = \frac{\det(\m{L}_Y)}{\det(\m{L} + \m{I})}, \quad Y \subseteq \mc{Y},
	\end{equation}
	where $\m{L}_Y := [\m{L}_{ij}]_{i,j \in Y}$ is the principal submatrix indexed by $Y$, $\m{I} \in \Rn{p \times p}$ is the identity matrix, and $\det(\m{L}_\emptyset) = 1$ by convention.
\end{mydef}
The normalizing constant in~\eqref{eq:dpp_L} follows from $\sum_{Y' \subseteq \mc{Y}} \det(\m{L}_{Y'}) = \det(\m{L} + \m{I})$; please refer to~\cite[Theorem 2.1]{Kulesza2012} for the proof.
The L-ensemble encodes the likelihood of each subset $Y$ through the determinant of the corresponding submatrix. Matrix $\m{L} \succeq \m{0}$ admits a Gram factorization $\m{L} = \m{G}^\top \m{G}$ for some $\m{G} \in \Rn{d \times p}$ with columns $\m{g}_1, \ldots, \m{g}_p$, where $\det(\m{L}_Y)$ equals the squared volume of the parallelepiped spanned by $\{\m{g}_i\}_{i \in Y}$~\cite{Tao2009}. Near orthogonal columns span a large parallelepiped, yielding high probability, while nearly collinear columns result in smaller $\mc{P}_{\m{L}}(Y)$. The matrix $\m{L}$ captures two attributes of each node. The first is \emph{quality}, measured by diagonal entries $\m{L}_{ii} = \|\m{g}_i\|^2$, which quantify the strength of node $i$. The second is \emph{diversity}, measured by entries $\m{L}_{ij} = \m{g}_i^\top \m{g}_j$, which quantify the similarity between nodes $i$ and $j$. Subsets that are simultaneously high quality and mutually diverse receive the highest probability under $\mc{P}_{\m{L}}$. The normalizing constant $\det(\m{L} + \m{I})$ aggregates the squared volumes of all $2^p$ possible subsets, ensuring that $\mc{P}_{\m{L}}$ is a valid probability measure on the power set $2^{\mc{Y}}$. 
\begin{mydef}[\textit{Marginal kernel~\cite{Kulesza2012}}]\label{def:marginal_kernel}
	The marginal kernel of an L-ensemble DPP is defined as
	\begin{equation}\label{eq:marginal_K}
		\m{K} := \m{L}(\m{L} + \m{I})^{-1} = \m{I} - (\m{L} + \m{I})^{-1}.
	\end{equation}
\end{mydef}

For any $A \subseteq \mc{Y}$, we have $\prob(A \subseteq Y) = \det(\m{K}_A)$. If $\lambda_1, \ldots, \lambda_p$ are the eigenvalues of $\m{L}$, then $\m{K}$ has eigenvalues $\kappa_i = \lambda_i / (\lambda_i + 1) \in [0,1)$. The marginal kernel \(\m{K}\) gives the node inclusion probabilities induced by the DPP. In particular, $\m{K}_{ii} = \prob(i \in Y)$ is the marginal probability that item $i$ is included in a random subset drawn from $\mc{P}_{\m{L}}$. Nodes with large kernel values $\m{L}_{ii}$ (high quality) have $\m{K}_{ii}$ close to one, while nodes with small values have $\m{K}_{ii}$ close to zero. Furthermore, $\m{K}_{ij}$ determines how likely nodes $i$ and $j$ are to be included together in a subset. The joint inclusion probability of nodes $i$ and $j$ satisfies $\prob(i, j \in Y) = \m{K}_{ii}\m{K}_{jj} - \m{K}_{ij}^2$, such that $\m{K}_{ij}^2$ quantifies the repulsion between nodes $i$ and $j$ (large when $\m{g}_i$ and $\m{g}_j$ are nearly collinear and vanishing when they are orthogonal). The eigenvalue mapping $\kappa_i = \lambda_i / (\lambda_i + 1)$ compresses the spectrum of $\m{L}$ into the interval $[0,1)$. In the DPP sampling algorithm~\cite{Hough2006,Kulesza2012}, the $i$-th eigenvector of $\m{L}$ can be selected independently with probability $\kappa_i$, where the expected cardinality of the selected subset is given as $\mbb{E}[|Y|] = \trace(\m{K}) = \sum_{i=1}^{p} \kappa_i$, which counts the number of effectively activated modes.

Furthermore, DPPs satisfy a negative dependence relation where for any $S_1, S_2 \subseteq \mc{Y}$ with $S_1 \cap S_2 = \emptyset$, the probability of including both $S_1$ and $S_2$ is upper bounded by the product of their individual inclusion probabilities as follows
\begin{equation}\label{eq:neg_dep}
	\prob(S_1 \cup S_2 \subseteq Y) \leq \prob(S_1 \subseteq Y) \cdot \prob(S_2 \subseteq Y).
\end{equation}
Inequality~\eqref{eq:neg_dep} follows from the determinantal structure of the marginal kernel given by $\prob(A \subseteq Y) = \det(\m{K}_A)$ satisfying $\det(\m{K}_{S_1 \cup S_2}) \leq \det(\m{K}_{S_1}) \det(\m{K}_{S_2})$ whenever $S_1 \cap S_2 = \emptyset$~\cite{Lyons2003}. Under independent selection the inequality holds with equality. 
Negative dependence is therefore the probabilistic mechanism through which the DPP enforces diversity, whereby each additional item included in $Y$ reduces the probability of including similar items. In \secref{sec:main_results}, we exploit this determinantal structure by constructing the L-ensemble matrix directly from the observability Gramian, such that the repulsive structure of the DPP aligns with the geometric redundancy among nodes. 

\section{From Gramians to Point Processes}\label{sec:main_results}
In this section, we establish that the observability Gramian $\m{W}_o$ is a valid L-ensemble DPP matrix. We then derive the marginal kernel $\m{K} = \m{W}_o(\m{W}_o + \m{I})^{-1}$, where $\m{K}_{ii}$ gives the inclusion probability of node $i$ and the eigenvalues of $\m{K}$ characterize modal activation. 
The corresponding controllability results follow from the duality stated in \secref{subsec:control}.

\subsection{Gramians and DPPs}\label{subsec:gramian_dpp}

We consider systems in which each node represents a state variable of the LTI system~\eqref{eq:model_DT}, that is, $\m{C}=\m{I}$. Let $\m{o}_i \in \Rn{Np}$ denote the $i$-th column of the observability matrix $\m{O}$~\eqref{eq:Obs_matrix}, such that the entry $(\m{W}_o)_{ij} = \m{o}_i^\top \m{o}_j$ quantifies the inner product of the observability contributions of states $i$ and $j$ over the horizon $N$, and let $\mr{Vol}(\{\m{o}_i\}_{i \in Y})$ denote the volume of the parallelepiped spanned by the columns indexed by $Y$. The following theorem shows that $\m{W}_o$ itself serves as a valid L-ensemble DPP matrix, and that the DPP probability of any node subset $Y$ is proportional to $\det((\m{W}_o)_{Y})$, the determinant of the corresponding principal submatrix.

\begin{theorem}[\textit{Gramian as a DPP}]\label{thm:gramian_dpp}
	Let $\m{L} := \m{W}_o \in \Rn{n \times n}$ with $\mc{Y}=\{1,\ldots,n\}$. Then $\m{L} \succeq \m{0}$, and the L-ensemble $\mc{P}_{\m{L}}$ is a valid probability measure on $2^{\mc{Y}}$ that satisfies, for every $Y \subseteq \mc{Y}$, the following
	\begin{equation}\label{eq:dpp_volume}
		\mc{P}_{\m{L}}(Y) = \frac{\det((\m{W}_o)_Y)}{\det(\m{W}_o + \m{I})} = \frac{\mr{Vol}^2\big(\{\m{o}_i\}_{i \in Y}\big)}{\det(\m{W}_o + \m{I})}.
	\end{equation}
	Moreover, $\mc{P}_{\m{L}}(Y) > 0$ if and only if the columns $\{\m{o}_i\}_{i \in Y}$ are linearly independent, and the system~\eqref{eq:model_DT} is observable over the horizon $N$ if and only if $\mc{P}_{\m{L}}(Y) > 0$ for every $Y \subseteq \mc{Y}$.
\end{theorem}
\begin{proof}
	Since $\m{W}_o = \m{O}^\top \m{O}$, positive semidefiniteness $\m{L} \succeq \m{0}$ follows immediately. Let $\m{E}_Y := [\m{e}_i]_{i \in Y} \in \Rn{n \times |Y|}$ collect the standard basis vectors $\m{e}_i \in \Rn{n}$ indexed by $Y$. Then $\m{L}_Y = (\m{O}\m{E}_Y)^\top (\m{O}\m{E}_Y)$ is the Gram matrix of the columns $\{\m{o}_i\}_{i \in Y}$, and the Gram determinant identity gives $\det(\m{L}_Y) = \mr{Vol}^2(\{\m{o}_i\}_{i \in Y})$~\cite[Sec.~2.2]{Kulesza2012}. For the normalizing identity, the $i$-th column of $\m{L} + \m{I}$ is $\m{e}_i + \m{l}_i$, where $\m{l}_i$ is the $i$-th column of $\m{L}$. Expanding $\det(\m{L}+\m{I})$ by multilinearity of the determinant in each column produces one term per subset $Y' \subseteq \mc{Y}$, and cofactor expansion of each term along its unit columns $\m{e}_j$, $j \notin Y'$, yields $\det(\m{L} + \m{I}) = \sum_{Y' \subseteq \mc{Y}} \det(\m{L}_{Y'})$, where the sum ranges over all $2^n$ subsets and $\det(\m{L}_\emptyset) = 1$~\cite[Theorem~2.1]{Kulesza2012}. Therefore, $\mc{P}_{\m{L}}$ sums to one over $2^{\mc{Y}}$, and~\eqref{eq:dpp_volume} holds. The Gram matrix $\m{L}_Y \succeq \m{0}$ is nonsingular if and only if the vectors $\{\m{o}_i\}_{i \in Y}$ are linearly independent~\cite[Thm.~7.2.10]{Horn1985}; it follows that $\det(\m{L}_Y) > 0$, and hence $\mc{P}_{\m{L}}(Y) > 0$, if and only if this independence holds. Furthermore, the system is observable over the horizon $N$ if and only if $\rank(\m{O}) = n$, equivalently the $n$ columns of $\m{O}$ are linearly independent. For this case every configuration $Y \subseteq \mc{Y}$ collects linearly independent columns and satisfies $\mc{P}_{\m{L}}(Y) > 0$. Conversely, $\mc{P}_{\m{L}}(\mc{Y}) > 0$ gives $\det(\m{W}_o) > 0$, and thus $\rank(\m{O}) = n$. This concludes the proof.
\end{proof}

\theoref{thm:gramian_dpp} establishes that the observability Gramian directly defines a valid DPP over the node index set $\mc{Y}=\{1,\ldots,n\}$. This connection provides a probabilistic interpretation of node selection. The DPP measure $\mc{P}_{\m{L}}$ assigns higher probability to node subsets $Y$ for which the principal submatrix $\m{L}_Y$ has a large determinant, equivalently, to configurations that are both observable and diverse. From \theoref{thm:gramian_dpp}, the rank condition admits a probabilistic characterization, whereby the system is observable over the horizon $N$ if and only if every node configuration receives positive probability under $\mc{P}_{\m{L}}$. We note that Equation~\eqref{eq:dpp_volume} does not introduce randomness in the system dynamics; it assigns a deterministic nonnegative weight $\det(\m{L}_Y)$ to each subset $Y$. The normalization in~\eqref{eq:dpp_L}, equivalently $\det(\m{L}+\m{I})=\sum_{Y'\subseteq\mc{Y}}\det(\m{L}_{Y'})$, converts these weights into a probability mass function over all possible subsets of the ground set $\mc{Y}$. For a linear measurement model in which a sensor measures one or multiple state variables across different nodes, that is, $\m{C} \neq \m{I}$, the selection is over the $p$ candidate sensors, and its objective $\det(\m{W}_o(Y))$ is not a principal minor of $\m{W}_o$. We therefore establish the corresponding DPP over the measurements of the observability matrix in~\appref{app:generalC}.

\begin{exmpl}[Illustrative example]\label{rem:example_dpp}
Consider an LTI system with $n=3$ states, and observability Gramian given by $\m{L}=\m{W}_o = \begin{bsmallmatrix} 4 & 2 & 0 \\ 2 & 3 & 0 \\ 0 & 0 & 1 \end{bsmallmatrix}$. By \theoref{thm:gramian_dpp}, the normalizing constant is $\det(\m{L}+\m{I})=32$, and every subset $Y\subseteq\{1,2,3\}$ receives probability $\mc{P}_{\m{L}}(Y)=\det(\m{L}_Y)/32$. Computing all $2^3$ subsets yields
\[
\begingroup
\setlength{\tabcolsep}{2.6pt}
\renewcommand{\arraystretch}{1.05}
\resizebox{\columnwidth}{!}{$
\begin{array}{c|cccccccc}
Y & \emptyset & \{1\} & \{2\} & \{3\} & \{1,2\} & \{1,3\} & \{2,3\} & \{1,2,3\} \\
\hline
\det(\m{L}_Y) & 1 & 4 & 3 & 1 & 8 & 4 & 3 & 8 \\
\mc{P}_{\m{L}}(Y) & \tfrac{1}{32} & \tfrac{4}{32} & \tfrac{3}{32} & \tfrac{1}{32} & \tfrac{8}{32} & \tfrac{4}{32} & \tfrac{3}{32} & \tfrac{8}{32}
\end{array}
$}
\endgroup
\]
which sums to one. Node $1$ has the highest singleton probability due to its largest diagonal entry $\m{L}_{11}=4$. Although nodes $1$ and $3$ are orthogonal ($\m{L}_{13}=0$), the pair $\{1,2\}$ is more likely than $\{1,3\}$ since the observability energy of nodes $1$ and $2$ outweighs the correlation penalty from $\m{L}_{12}=2$.
\end{exmpl}

\begin{figure}[t]
\centering
\resizebox{\linewidth}{!}{%
\begin{tikzpicture}[transform shape, every node/.style={font=\normalsize}]
			\begin{scope}
				\fill[LCSS!5]
					(0,0) -- (1.50,0.14) -- (1.67,1.62) -- (0.17,1.48) -- cycle;
				\draw[LCSS!45, densely dashed, thin]
					(1.50,0.14) -- (1.67,1.62) -- (0.17,1.48);
			\draw[-{Stealth[length=4pt]}, semithick, black!22] (0,0) -- (1.40,0.39)
				node[right, font=\scriptsize, black!35] {$\m{o}_2$};
			\draw[-{Stealth[length=4pt]}, semithick, black!22] (0,0) -- (0.39,1.18)
				node[above, font=\scriptsize, black!35] {$\m{o}_4$};
				\draw[-{Stealth[length=4.5pt]}, very thick, LCSS!85!black] (0,0) -- (1.50,0.14)
					node[right, font=\scriptsize, LCSS!85!black] {$\m{o}_1$};
				\draw[-{Stealth[length=4.5pt]}, very thick, LCSS!60!black] (0,0) -- (0.17,1.48)
					node[above, font=\scriptsize, LCSS!60!black] {$\m{o}_3$};
				\node[font=\scriptsize, LCSS!80!black] at (0.76,0.64) {$\mr{Vol}^2$};
			\fill[black] (0,0) circle (1.5pt);
				\node[font=\scriptsize] at (0.98,-0.45)
					{$Y_1 = \{1,3\}$: \textcolor{LCSS}{diverse}};
			\node[font=\scriptsize\bfseries, anchor=south] at (0.98,1.8)
				{(a) High $\det(\m{L}_{Y_1})$};
		\end{scope}
		\begin{scope}[shift={(4,0)}]
				\begin{scope}
					\clip (-0.1,-0.1) rectangle (2.9,2.1);
						\fill[red!5]
							(0,0) -- (1.50,0.14) -- (2.73,0.66) -- (1.23,0.52) -- cycle;
						\draw[red!38, densely dashed, thin]
							(1.50,0.14) -- (2.73,0.66) -- (1.23,0.52);
				\end{scope}
				\draw[-{Stealth[length=4pt]}, semithick, black!22] (0,0) -- (0.17,1.48)
					node[above, font=\scriptsize, black!35] {$\m{o}_3$};
				\draw[-{Stealth[length=4pt]}, semithick, black!22] (0,0) -- (0.39,1.18)
					node[above right, font=\scriptsize, black!35] {$\m{o}_4$};
					\draw[-{Stealth[length=4.5pt]}, very thick, red!80!black] (0,0) -- (1.50,0.14)
						node[below right, font=\scriptsize, red!80!black] {$\m{o}_1$};
					\draw[-{Stealth[length=4.5pt]}, very thick, red!72!black] (0,0) -- (1.23,0.52)
						node[above left, font=\scriptsize, red!72!black] {$\m{o}_2$};
					\node[font=\scriptsize, red!72!black, rotate=14] at (1.44,0.34)
						{$\mr{Vol}^2$};
				\fill[black] (0,0) circle (1.5pt);
				\node[font=\scriptsize] at (0.98,-0.45)
					{$Y_2 = \{1,2\}$: \textcolor{red!72!black}{redundant}};
				\node[font=\scriptsize\bfseries, anchor=south] at (0.98,1.8)
					{(b) Low $\det(\m{L}_{Y_2})$};
			\end{scope}
			\node[font=\scriptsize\itshape] at (6.0,1.0)
				{$\mc{P}_{\m{L}}(Y)\propto\det(\m{L}_Y)$};
	\end{tikzpicture}%
	}
	\vspace{-0.3cm}
\caption{Geometric view of the Gramian DPP equivalence (\theoref{thm:gramian_dpp}). The vectors $\m{o}_i$ denote the columns of $\m{O}$; bold arrows are the selected subset $Y$. (a)~Near orthogonal pair $\{1,3\}$ yields a large $\det((\m{W}_o)_{Y_1})$, which produces a high $\mc{P}_{\m{L}}(Y_1)$. (b)~Nearly collinear pair $\{1,2\}$ collapses the determinant, and thus suppresses $\mc{P}_{\m{L}}(Y_2)$.}
\label{fig:dpp_geometry}
\end{figure}

We note here that $\det(\m{L}_Y)$, which defines the subset likelihood in the DPP, coincides with the objective of \probref{prob:node_selection}. %

From \theoref{thm:gramian_dpp}, the Gram structure of $\m{L}$ offers a geometric interpretation of probability $\mc{P}_{\m{L}}(Y)$ as the normalized squared volume of the parallelepiped spanned by the columns of $\m{O}$ indexed by $Y$. Let $q_i := \sqrt{(\m{W}_o)_{ii}} \geq 0$ denote the individual observability energy of node $i$, and define the normalized entries $\m{\phi}_{ij} := (\m{W}_o)_{ij} / (q_i q_j)$ for all $i,j \in Y$ with $q_i,q_j>0$. Then, $\mc{P}_{\m{L}}(Y)$ is proportional to
	\begin{equation}\label{eq:quality_diversity}
		\mc{P}_{\m{L}}(Y) \propto \left(\prod_{i \in Y} q_i^2\right) \det(\m{S}_Y),
	\end{equation}
	where $\m{S}_Y := [\m{\phi}_{ij}]_{i,j \in Y}$ is the normalized observability Gramian submatrix of the configuration $Y$. This decomposition~\eqref{eq:quality_diversity} reveals two competing objectives in node selection, that is, selecting nodes with large individual Gramian energy $q_i^2 = (\m{W}_o)_{ii}$, while ensuring that their observability directions are mutually orthogonal; see \figref{fig:dpp_geometry}. The quality term $q_i^2$ measures how much observability energy node $i$ contributes in isolation, while the diversity term $\det(\m{S}_Y)$ equals one when the selected observability vectors are orthogonal and approaches zero when they become collinear, thereby penalizing configurations where multiple nodes observe similar linear combinations of the state variables. This decomposition motivates the following definition.

\begin{mydef}[\textit{Diversity}]\label{def:diversity}
	For a node subset $Y \subseteq \mc{Y}$ with $q_i = \sqrt{(\m{W}_o)_{ii}} > 0$ for all $i \in Y$, the \emph{diversity index} of the configuration $Y$ is defined as
	\begin{equation}\label{eq:diversity_index}
		\mc{D}(Y) := \det(\m{S}_Y), \quad \mc{D}(Y) \in [0,1],
	\end{equation}
		where $\m{S}_Y := [(\m{W}_o)_{ij}/(q_i\, q_j)]_{i,j \in Y}$ is the normalized observability Gramian of the configuration $Y$.
\end{mydef}

 \defref{def:diversity} formalizes the notion of diversity, where $\mc{D}(Y) = 1$ when the selected nodes induce mutually orthogonal observability directions, and $\mc{D}(Y) = 0$ when these directions are linearly dependent. A sensor configuration $Y$ is diverse when its nodes contribute linearly independent measurements of the state, and is redundant when their observability directions overlap. This geometric interpretation motivates a spectral characterization of subset inclusion through the marginal kernel, which is introduced in the following subsection.

\parstartc{Revisiting~\exmpref{rem:example_dpp}} The individual observability energies are $q_1^2 = 4$, $q_2^2 = 3$, $q_3^2 = 1$. For the pair $\{1,2\}$, the quality term equals $q_1^2 q_2^2 = 12$ and the diversity index is $\mc{D}(\{1,2\}) = 2/3$, yielding $\det(\m{L}_{\{1,2\}}) = 8$. In contrast, the pair $\{1,3\}$ induces mutually orthogonal observability directions and attains $\mc{D}(\{1,3\}) = 1$, while the product of its observability energies is $q_1^2 q_3^2 = 4$, and thus $\det(\m{L}_{\{1,3\}}) = 4$. The DPP assigns a higher probability to $\{1,2\}$, since the larger observability energies outweigh the reduced diversity index. The decomposition~\eqref{eq:quality_diversity} thus balances the two competing objectives.

The following proposition relates the objective of \probref{prob:node_selection} to the accuracy with which the state variables are estimated from noisy measurements.

\begin{myprs}[\textit{Estimation error covariance}]\label{prop:estimation}
	Let $\m{W}_o \succ \m{0}$ and consider the measurements $\m{y}[k] = \m{x}[k] + \m{v}[k]$ over the horizon $N$, where $\m{v}[k]$ is zero mean Gaussian with covariance $\sigma^2 \m{I}$. Let $\widehat{\m{x}}[0]$ be the least squares estimate of $\m{x}[0]$, define the estimation error $\widetilde{\m{x}}[0] := \widehat{\m{x}}[0]-\m{x}[0]$, and let $\m{\Sigma} := \sigma^2 \m{W}_o^{-1}$ denote its covariance. For every configuration $Y \subseteq \mc{Y}$ with complement $\bar{Y} = \mc{Y} \setminus Y$,
	\begin{equation}\label{eq:jacobi}
		\det\big((\m{W}_o)_Y\big) = \det(\m{W}_o)\, \det\big((\m{W}_o^{-1})_{\bar{Y}}\big),
	\end{equation}
	and the error covariance of $\widetilde{\m{x}}_Y[0]$ conditioned on $\widetilde{\m{x}}_{\bar{Y}}[0]$ is
	\begin{equation}\label{eq:cond_cov}
		\mr{cov}\big(\widetilde{\m{x}}_Y[0] \mid \widetilde{\m{x}}_{\bar{Y}}[0]\big) = \sigma^2 \big((\m{W}_o)_Y\big)^{-1}.
	\end{equation}
\end{myprs}

\begin{proof}
	The proof of~\eqref{eq:jacobi} follows from the identity of Jacobi for the minors of an inverse matrix~\cite[Sec.~0.8.4]{Horn1985}, which gives $\det((\m{W}_o^{-1})_{\bar{Y}}) = \det((\m{W}_o)_Y)/\det(\m{W}_o)$ for every $Y \subseteq \mc{Y}$. For~\eqref{eq:cond_cov}, stacking the measurements over the horizon $N$ yields $\m{Y} = \m{O}\m{x}[0] + \bm{\eta}$, where $\bm{\eta}$ collects the noise over the horizon and has covariance $\sigma^2\m{I}$. The least squares estimation error $\widetilde{\m{x}}[0]$ is zero mean Gaussian with covariance $\m{\Sigma} = \sigma^2 (\m{O}^\top \m{O})^{-1} = \sigma^2 \m{W}_o^{-1}$. Partitioning $\widetilde{\m{x}}[0]$ according to $Y$ and $\bar{Y}$, the conditional error covariance is the Schur complement of $\m{\Sigma}_{\bar{Y}}$ in $\m{\Sigma}$, given by $\m{\Sigma}_{Y} - \m{\Sigma}_{Y\bar{Y}} \m{\Sigma}_{\bar{Y}}^{-1} \m{\Sigma}_{\bar{Y}Y}$. From the block inversion of $\m{\Sigma}$ with respect to $Y$ and $\bar{Y}$~\cite[Sec.~0.7.3]{Horn1985}, this Schur complement equals $((\m{\Sigma}^{-1})_Y)^{-1} = \sigma^2 ((\m{W}_o)_Y)^{-1}$, which yields~\eqref{eq:cond_cov}. This concludes the proof.
\end{proof}

Under the measurement noise considered in \propref{prop:estimation}, the matrix $\m{W}_o/\sigma^2$ is the Fisher information matrix of $\m{x}[0]$, and the determinant of the error covariance quantifies the volume of the confidence ellipsoid of the estimate; determinant measures of this form are used for sensor selection in~\cite{Joshi2009,Liu2023}. From~\eqref{eq:jacobi}, and given that $\det(\m{W}_o)$ does not depend on $Y$. Then maximizing $\logdet{(\m{W}_o)_Y}$ is equivalent to maximizing $\det(\m{\Sigma}_{\bar{Y}})$, the determinant of the error covariance of the remaining state variables. This equivalence involves the marginal covariance $\m{\Sigma}_{\bar{Y}}$ and assumes no knowledge of the remaining state variables. As such, the selected configuration carries the most information about its own state variables, and the remaining ones retain the largest share of the estimation uncertainty. From~\eqref{eq:cond_cov}, the same objective minimizes the determinant of the conditional error covariance of the selected state variables, for which the confidence ellipsoid of the estimate results in the smallest volume. When $\m{W}_o$ is diagonal, maximizing $\logdet{(\m{W}_o)_Y}$ and $\det(\m{\Sigma}_{\bar{Y}})$ reduce to selecting the $r$ nodes of largest energy $(\m{W}_o)_{ii}$. Note that the covariance in~\eqref{eq:cond_cov} is the Schur complement of the error covariance with respect to the remaining state variables. This is also the case for the DPP where the inclusion probability of a node given an already selected configuration is the Schur complement of the marginal kernel (\coref{cor:cond_incl}).

\subsection{The marginal kernel as a modal observability filter}\label{subsec:marginal}
The eigenvectors $\m{v}_1, \ldots, \m{v}_n$, defined as the columns of $\m{V}$ in the eigendecomposition $\m{L}=\m{V}\m{\Lambda}\m{V}^\top$, define the \emph{observability modes} of the system. The corresponding eigenvalues $\lambda_1 \geq \cdots \geq \lambda_n \geq 0$ measure the energy concentrated in each mode. The marginal kernel $\m{K} = \m{L}(\m{L} + \m{I})^{-1}$ preserves these modal directions while mapping each eigenvalue through $\lambda_i \mapsto \kappa_i = \lambda_i/(\lambda_i + 1) \in [0,1)$; it converts unbounded Gramian energy into bounded activation probabilities accordingly.

\begin{theorem}[\textit{Marginal kernel and modal activation}]\label{thm:marginal_modal}
	Let $\m{L}= \m{W}_o$ be the DPP matrix with eigendecomposition $\m{L} = \m{V} \m{\Lambda} \m{V}^\top$, where $\m{V} \in \Rn{n \times n}$ is orthogonal, $\m{\Lambda} = \mr{diag}(\lambda_1, \ldots, \lambda_n)$, and $\lambda_1 \geq \cdots \geq \lambda_n \geq 0$. Then, the marginal kernel written as
	\begin{equation}\label{eq:K_modal}
		\m{K} = \m{L}(\m{L} + \m{I})^{-1} = \m{V}\, \mr{diag}\!\left(\frac{\lambda_1}{\lambda_1 + 1}, \ldots, \frac{\lambda_n}{\lambda_n + 1}\right) \m{V}^\top,
	\end{equation}
	has the following properties.
	\begin{enumerate}
		\item[\emph{(i)}] The diagonal entry $\m{K}_{ii} = \prob(i \in Y)$ is the marginal probability that node $i$ is included in a random subset drawn from $\mc{P}_{\m{L}}$. 
		\item[\emph{(ii)}] The eigenvalues $\kappa_i \hspace{-0.1cm}=\hspace{-0.1cm} \lambda_i / (\lambda_i + 1)\hspace{-0.1cm} \in \hspace{-0.1cm}[0, 1)$ are the inclusion probabilities of the $i$-th \emph{observability mode}, where mode $i$ is selected independently with probability $\kappa_i$.
		\item[\emph{(iii)}] The expected number of selected nodes is given by
		\begin{equation}\label{eq:expected_card}
			\mbb{E}[|Y|] = \trace(\m{K}) = \sum_{i=1}^{n} \frac{\lambda_i}{\lambda_i + 1},
		\end{equation}
		and the variance is $\mr{Var}[|Y|] = \sum_{i=1}^{n} \lambda_i / (\lambda_i + 1)^2$.
	\end{enumerate}
\end{theorem}

\begin{proof}
	With $\m{L} \succeq \m{0}$, this orthogonal diagonalization holds even when eigenvalues are repeated, and given that $\m{L} = \m{V}\m{\Lambda}\m{V}^\top$, applying the scalar function $f(\lambda) = \lambda/(\lambda+1)$ to the spectrum gives $\m{K} = \m{V}\,f(\m{\Lambda})\,\m{V}^\top$, which yields the spectral form~\eqref{eq:K_modal}. For \textit{(i)}, by \defref{def:marginal_kernel}, $\prob(i \in Y) = \det(\m{K}_{\{i\}}) = \m{K}_{ii}$. For \textit{(ii)}, the DPP sampling algorithm~\cite{Hough2006,Kulesza2012} selects eigenvector $\m{v}_i$ independently with probability $\kappa_i = \lambda_i/(\lambda_i + 1)$ via a Bernoulli trial, and then projects the selected eigenvectors onto node coordinates. For \textit{(iii)}, linearity of expectation gives $\mbb{E}[|Y|] = \sum_{i=1}^{n} \prob(i \in Y) = \trace(\m{K})$, and since the mode selections are independent Bernoulli variables, $\mr{Var}[|Y|] = \sum_{i=1}^{n} \kappa_i(1 - \kappa_i) = \sum_{i=1}^{n} \lambda_i/(\lambda_i+1)^2$. This concludes the proof.
\end{proof}

From the above theorem, the expected cardinality~\eqref{eq:expected_card} defines an \emph{effective observable rank} as $n_{\mr{eff}} := \mbb{E}[|Y|] = \trace(\m{K}) = \sum_{i=1}^{n} \frac{\lambda_i}{\lambda_i + 1}$, which counts the number of observability modes that are significantly activated. Modes with large eigenvalues contribute $\kappa_i \approx 1$, while modes with small eigenvalues contribute $\kappa_i \approx 0$. As such, $n_{\mr{eff}}$ satisfies $0 \leq n_{\mr{eff}} \leq \rank(\m{L})$ and quantifies the effective dimension of the observable subspace. The value $n_{\mr{eff}}$ characterizes the expected number of nodes selected by the DPP. We note that the mapping $\lambda_i \mapsto \lambda_i/(\lambda_i+1)$ normalizes the modal observability energies to $\kappa_i \in [0,1)$, where the $\m{I}$ term in $\m{K} = \m{W}_o(\m{W}_o+\m{I})^{-1}$ introduces a unit reference. As such, a mode is activated when its observability energy $\lambda_i$ exceeds this reference.

\parstartc{Revisiting~\exmpref{rem:example_dpp}} The marginal kernel diagonal entries are $\m{K}_{11} = 3/4$, $\m{K}_{22} = 11/16$, and $\m{K}_{33} = 1/2$, and thus the DPP selects $n_{\mr{eff}} = \trace(\m{K}) = 31/16 \approx 1.94$ nodes in expectation, whereas $\rank(\m{W}_o) = 3$. The correlation between nodes $1$ and $2$ lowers the effective observable rank. The same observability energies with $(\m{W}_o)_{12} = 0$ yield $n_{\mr{eff}} = 2.05$. Node~$1$ attains the highest inclusion probability, consistent with its largest observability energy $(\m{W}_o)_{11} = 4$.

The following corollary complements the above results by quantifying how the inclusion probability of a single node changes once a set of nodes has already been selected.

\begin{mycor}[\textit{Conditional inclusion}]\label{cor:cond_incl}
	Let $S \subset \mc{Y}$ be a set of already selected nodes with $\det(\m{K}_S) > 0$, and let $j \notin S$. The conditional inclusion probability is given by
	\begin{equation}\label{eq:cond_incl}
		\prob(j \in Y \mid S \subseteq Y) = \m{K}_{jj} - \m{K}_{j,S}\, \m{K}_S^{-1}\, \m{K}_{S,j},
	\end{equation}
	where $\m{K}_{j,S} \in \Rn{1 \times |S|}$ and $\m{K}_S \in \Rn{|S| \times |S|}$ are the corresponding submatrices of $\m{K}$.
\end{mycor}

\begin{proof}
	By the inclusion property (\defref{def:marginal_kernel}), $\prob(S \cup \{j\} \subseteq Y) = \det(\m{K}_{S \cup \{j\}})$ and $\prob(S \subseteq Y) = \det(\m{K}_S)$. By Bayes' rule, $\prob(j \in Y \mid S \subseteq Y) = \det(\m{K}_{S \cup \{j\}}) / \det(\m{K}_S)$, which equals the Schur complement of $\m{K}_{jj}$ in $\m{K}_{S \cup \{j\}}$ with respect to $\m{K}_S$~\cite{Kulesza2012}. This yields~\eqref{eq:cond_incl} and concludes the proof.
\end{proof}

Equation~\eqref{eq:cond_incl} quantifies how the inclusion probability of node $j$ decreases as more nodes in $S$ are selected. The term $\m{K}_{j,S}\, \m{K}_S^{-1}\, \m{K}_{S,j}$ measures the observability relation of node $j$ to the already selected set $S$ through the marginal kernel. When the Gramian row of node $j$ is nearly a linear combination of the rows indexed by $S$, the conditional probability drops. This is relevant in practice when additional sensor nodes are to be selected after an initial set $S$ has been chosen.

\parstartc{Revisiting~\exmpref{rem:example_dpp}} $\prob(2 \in Y \mid 1 \in Y) = \m{K}_{22} - \m{K}_{21}\m{K}_{11}^{-1}\m{K}_{12} = 2/3$, which is less than the marginal $\prob(2 \in Y) = 11/16$; including node~$1$ reduces the benefit of adding the correlated node~$2$. In contrast, $\prob(3 \in Y \mid 1 \in Y) = \m{K}_{33} = 1/2 = \prob(3 \in Y)$, since $\m{K}_{13} = 0$ and nodes~$1$ and~$3$ have orthogonal Gramian rows.

\section{Inclusion Properties and Node Selection}\label{sec:node_selection}
Having established the Gramian DPP relation and the marginal kernel, we now derive properties of the DPP that have direct implications for sensor node selection. 
We then recover the submodularity of $\logdett$ from a geometric argument and reformulate \probref{prob:node_selection} as a MAP problem.

\subsection{Inclusion monotonicity and node redundancy}\label{subsec:structural}
The first result establishes that if one system is more observable than another in the positive semidefinite sense, then the DPP assigns a higher inclusion probability to every sensor configuration under the more observable system.

\begin{mycor}[\textit{Inclusion monotonicity}]\label{cor:stoch_dom}
	Let $\m{L}^{(1)} = \m{W}_o^{(1)}$ and $\m{L}^{(2)} = \m{W}_o^{(2)}$. If $\m{W}_o^{(1)} \preceq \m{W}_o^{(2)}$, then for every sensor configuration $S \subseteq \mc{Y}$, the inclusion probability under the more observable system is at least as large, that is,
	\begin{equation}\label{eq:stoch_dom}
		\prob_{\m{L}^{(1)}}(S \subseteq Y) \leq \prob_{\m{L}^{(2)}}(S \subseteq Y), \quad \forall\, S \subseteq \mc{Y},
	\end{equation}
	where $\prob_{\m{L}^{(i)}}(S \subseteq Y) = \det(\m{K}_S^{(i)})$ is the probability that the configuration $S$ is included in a random subset $Y$ drawn from the DPP with kernel $\m{L}^{(i)}$.
\end{mycor}
\begin{proof}
	The function $f(\lambda) = \lambda/(\lambda+1)$ is operator monotone on $[0,\infty)$. Therefore $\m{L}^{(1)} \preceq \m{L}^{(2)}$ implies
	\begin{equation*}
		\m{K}^{(1)} = f\!\left(\m{L}^{(1)}\right) \preceq f\!\left(\m{L}^{(2)}\right) = \m{K}^{(2)},
	\end{equation*}
	such that $\m{0} \preceq \m{K}^{(1)} \preceq \m{K}^{(2)} \preceq \m{I}$. For any $S \subseteq \mc{Y}$, the principal submatrix retains the ordering $\m{K}_S^{(1)} \preceq \m{K}_S^{(2)}$. If $\m{K}_S^{(1)}$ is singular, then $\det(\m{K}_S^{(1)}) = 0 \leq \det(\m{K}_S^{(2)})$. Otherwise, $\m{K}_S^{(1)}$ is nonsingular; then congruence by $(\m{K}_S^{(1)})^{-1/2}$ applied to $\m{K}_S^{(2)} - \m{K}_S^{(1)} \succeq \m{0}$ gives $\big(\m{K}_S^{(1)}\big)^{-1/2}\m{K}_S^{(2)}\big(\m{K}_S^{(1)}\big)^{-1/2} \succeq \m{I}$. Thus, $\det(\m{K}_S^{(2)}) \geq \det(\m{K}_S^{(1)})$. Now for both cases, $\prob_{\m{L}^{(1)}}(S \subseteq Y) = \det(\m{K}_S^{(1)}) \leq \det(\m{K}_S^{(2)}) = \prob_{\m{L}^{(2)}}(S \subseteq Y)$; this yields~\eqref{eq:stoch_dom} and concludes the proof.
\end{proof}

\coref{cor:stoch_dom} establishes an \textit{inclusion monotonicity} property for the Gramian DPP. Every inclusion event $\{S \subseteq Y\}$ has at least as high a probability under the larger kernel, so improving observability raises the inclusion probability of every configuration $S \subseteq \mc{Y}$. When the two kernels commute, this ordering extends to a stochastic dominance\footnote{A probability measure over subsets is stochastically dominant over another if it assigns at least as much probability to every increasing event, that is, every event that remains satisfied when elements are added to the subset~\cite{Lyons2003}.} of $\mc{P}_{\m{L}^{(2)}}$ over $\mc{P}_{\m{L}^{(1)}}$~\cite[Thm~8.1]{Lyons2003}. The ordering $\m{L}^{(1)} \preceq \m{L}^{(2)}$ holds whenever the Gramian of system~2 dominates that of system~1, which occurs whenever we modify $(\m{A},\m{C})$ to improve observability, thereby monotonically increasing the inclusion probabilities of all configurations. The next result introduces a direct interpretation of the unselected nodes as a DPP with kernel $\m{I} - \m{K}$, which characterizes what is unobserved by the selected configuration.

\begin{mycor}[\textit{Complement characterization}]\label{cor:hodge}
	Let $\m{L} = \m{W}_o$ and let $Y \sim \mc{P}_{\m{L}}$ be the random subset of selected sensor nodes. The complement $\bar{Y} := \mc{Y} \setminus Y$, representing the unselected nodes, is itself a DPP with marginal kernel $\m{I} - \m{K}$, where $\m{K} = \m{W}_o(\m{W}_o + \m{I})^{-1}$. That is, the complement inclusion probability satisfies
	\begin{equation}\label{eq:hodge}
		\prob(\bar{Y} \supseteq S) = \det\big((\m{I} - \m{K})_S\big), \quad \forall\, S \subseteq \mc{Y}.
	\end{equation}
	The diagonal entry $(\m{I}-\m{K})_{ii} = 1-\m{K}_{ii}$ is the probability that node $i$ is excluded from the selected set.
\end{mycor}

\begin{proof}
	Since $\m{K} = \m{W}_o(\m{W}_o + \m{I})^{-1}$, the eigenvalues of $\m{K}$ satisfy $\kappa_i = \lambda_i/(\lambda_i+1) \in [0,1)$, so $\m{I} - \m{K}$ is positive semidefinite with eigenvalues $1-\kappa_i = (\lambda_i+1)^{-1} \in (0,1]$. By the complement property of DPPs~\cite{Tao2009}, if $Y$ is a DPP with marginal kernel $\m{K}$, then $\bar{Y}$ is a DPP with marginal kernel $\m{I}-\m{K}$, yielding~\eqref{eq:hodge}.
\end{proof}

The above corollary provides a characterization of the unobserved portion of the state-space. At the node level, exclusion is quantified by $1-\m{K}_{ii}$. In the modal basis, the complement kernel has eigenvalues $1-\kappa_i = (\lambda_i+1)^{-1}$, so modes with large observability eigenvalues have small modal exclusion probabilities, while weak modes remain largely unobserved. The kernel $\m{I} - \m{K}$ therefore quantifies the residual observability that remains unobserved by the selected sensor configuration. This complements \propref{prop:estimation}, where the objective of \probref{prob:node_selection} is written through the error covariance of the state variables that are not selected.

\parstartc{Revisiting~\exmpref{rem:example_dpp}} The exclusion probabilities are $1-\m{K}_{11} = 1/4$, $1-\m{K}_{22} = 5/16$, and $1-\m{K}_{33} = 1/2$. Node~$3$ has the highest exclusion probability, and this confirms that weakly observable nodes dominate the complement DPP.

\begin{mycor}[\textit{Negative dependence}]\label{cor:neg_dep}
	Let $\m{L} = \m{W}_o$ and $Y \sim \mc{P}_{\m{L}}$. For any two disjoint sensor configurations $S_1, S_2 \subseteq \mc{Y}$, the joint inclusion probability satisfies
	\begin{equation}\label{eq:neg_corr}
		\prob(S_1 \cup S_2 \subseteq Y) \leq \prob(S_1 \subseteq Y) \cdot \prob(S_2 \subseteq Y).
	\end{equation}
	In particular, for any two nodes $i \neq j$,
	$
		\prob(i, j \in Y) = \m{K}_{ii}\m{K}_{jj} - \m{K}_{ij}^2 \leq \prob(i \in Y) \cdot \prob(j \in Y),
	$
	where the gap $\m{K}_{ij}^2$ increases with the observability correlation between nodes $i$ and $j$.
\end{mycor}

\begin{proof}
	By the inclusion property, $\prob(A \subseteq Y) = \det(\m{K}_A)$ for any $A \subseteq \mc{Y}$. For disjoint $S_1, S_2$, Fischer's inequality for positive semidefinite matrices gives $\det(\m{K}_{S_1 \cup S_2}) \leq \det(\m{K}_{S_1})\det(\m{K}_{S_2})$~\cite{Lyons2003}, which yields~\eqref{eq:neg_corr}. The pairwise identity follows from $\det(\m{K}_{\{i,j\}}) = \m{K}_{ii}\m{K}_{jj} - \m{K}_{ij}^2$.
\end{proof}

The aforementioned result implies that the DPP \emph{avoids redundancy}. Selecting node $i$ reduces the probability of also selecting node $j$ by $\m{K}_{ij}^2$. Nodes where $\m{K}_{ij}^2$ is close to $\m{K}_{ii}\m{K}_{jj}$ have a small joint inclusion probability and are strongly repelled. The DPP framework enables sampling while avoiding redundancy as compared to independent sampling that does not account for diversity.

\parstartc{Revisiting~\exmpref{rem:example_dpp}} $\prob(1,2 \in Y) = 1/2$ is strictly less than $\prob(1 \in Y)\prob(2 \in Y) = 33/64$, with a gap equal to $\m{K}_{12}^2 = 1/64$. In contrast, $\prob(1,3 \in Y) = 3/8 = \prob(1 \in Y)\prob(3 \in Y)$, so nodes~$1$ and~$3$ exhibit no repulsion because $\m{K}_{13} = 0$; orthogonal Gramian rows contribute independently to the DPP.

Based on the above results, a probabilistic interpretation of node selection is provided. The system dynamics and the Gramian remain deterministic; the probability model enters through the DPP measure in~\eqref{eq:dpp_L} and its Gramian form in~\eqref{eq:dpp_volume}, which define a distribution over subsets $Y \subseteq \mc{Y}$. This interpretation supports solving the combinatorial \probref{prob:node_selection} by using DPP properties to characterize optimal node configurations.

\subsection{Submodularity and MAP reformulation}\label{subsec:submod}
We now show that the $\logdett$ objective retains a diminishing returns property from the DPP matrix geometry. The following lemma establishes its submodularity and monotonicity.

\begin{mylem}[\textit{Submodularity and monotonicity}]\label{lem:submod}
	The set function $f(Y) := \log \det(\m{L}_Y)$ is submodular on $\mc{Y}$, with $f(\emptyset)=0$. That is, for any $A \subseteq B \subseteq \mc{Y}$ and $j \notin B$,
	\begin{equation}\label{eq:submod}
		f(A \cup \{j\}) - f(A) \geq f(B \cup \{j\}) - f(B).
	\end{equation}
	Moreover, if $\m{L} \succeq \m{I}$, then $f$ is monotone nondecreasing, that is, $f(A) \leq f(A \cup \{j\})$ for all $A \subseteq \mc{Y}$ and $j \notin A$.
\end{mylem}

\begin{proof}
	Let $\m{L}=\m{G}^\top \m{G}$ with columns $\{\m{g}_i\}_{i\in\mc{Y}}$. For $S\subseteq\mc{Y}$ and $j\notin S$, the marginal gain can be written as
	\begin{equation*}
		f(S\cup\{j\})-f(S)=\log h_j^2(S),
	\end{equation*}
	where $h_j^2(S):=(\m{L})_{jj}-\m{L}_{j,S}\m{L}_S^{-1}\m{L}_{S,j}$ is the Schur complement of $\m{L}_S$ in $\m{L}_{S\cup\{j\}}$~\cite{Kulesza2012}, with $h_j^2(\emptyset):=\m{L}_{jj}$. Equivalently,
	$ h_j^2(S)=\left\|(\m{I}-\Pi_S)\m{g}_j\right\|_2^2$,
	with $\Pi_S$ the orthogonal projector onto $\mr{span}\{\m{g}_i:i\in S\}$. If $A\subseteq B$, then $\mr{span}(A)\subseteq\mr{span}(B)$. Hence $\|(\m{I}-\Pi_B)\m{g}_j\|_2\leq\|(\m{I}-\Pi_A)\m{g}_j\|_2$, so $h_j(B)\leq h_j(A)$, and~\eqref{eq:submod} follows. For monotonicity, suppose $\m{L}\succeq\m{I}$ as stated in \lemref{lem:submod}. Then $\m{L}_{S\cup\{j\}}\succeq\m{I}$ and $\m{L}_{S\cup\{j\}}^{-1}\preceq\m{I}$. The block inverse formula gives $(\m{L}_{S\cup\{j\}}^{-1})_{jj}=h_j^{-2}(S)\leq1$, and thus $h_j^2(S)\geq1$. Therefore $\log h_j^2(S)\geq0$ for all $S\subseteq\mc{Y}$ and $j\notin S$. This concludes the proof.
\end{proof}

\lemref{lem:submod} recovers the submodularity of $\logdet{\m{L}_Y}$ from a geometric argument. The marginal gain of adding node $j$ equals the log Schur complement of $\m{L}_A$ in $\m{L}_{A\cup\{j\}}$, which decreases as more nodes are added to subset $A$; this exhibits the diminishing returns property that defines submodularity. We note that with $\m{C} = \m{I}$ the Gramian~\eqref{eq:LinObsGram} satisfies $\m{L} = \m{W}_o \succeq \m{I}$. The monotonicity condition of \lemref{lem:submod} is therefore satisfied, and the objective of \probref{prob:node_selection} is monotone submodular.

\begin{algorithm}[t]
	\caption{MAP greedy node selection for \probref{prob:node_selection}.}\label{alg:greedy_map}
	\DontPrintSemicolon
	\textbf{Input:} DPP kernel $\m{L} = \m{W}_o \in \Rn{n \times n}$, subset cardinality $r$\;\vspace{0.1cm}
	$Y \leftarrow \emptyset$\;
	\For{$t \leftarrow 1$ \KwTo $r$}{
		\textbf{select:} $j^\star \leftarrow \argmax_{j \notin Y}\; \m{L}_{jj} - \m{L}_{j,Y}\, \m{L}_Y^{-1}\, \m{L}_{Y,j}$
		// \textit{residual observability energy of node $j$ given $Y$}\;
		$Y \leftarrow Y \cup \{j^\star\}$\;
	}
	\textbf{Output:} node configuration $Y$ with $|Y| = r$\;
\end{algorithm}

Based on the results established above, the node selection \probref{prob:node_selection} can be reformulated as a $k$-DPP MAP inference problem. Restricting the process to subsets of fixed cardinality $|Y| = r$ yields the $k$-DPP, which can be sampled efficiently using elementary symmetric polynomial recursions~\cite{Kulesza2012}. The MAP problem of the $k$-DPP is equivalent to finding the most likely subset of size $r$. While this problem is NP-hard, greedy and local search heuristics are used in practice~\cite{Gillenwater2012}. The following proposition states the MAP reformulation of \probref{prob:node_selection}.
\begin{myprs}[\textit{$k$-DPP MAP reformulation}]\label{pro:kdpp_map}
For $\m{L} = \m{W}_o$ and a budget $r < n$, the $k$-DPP MAP problem is
\begin{equation}\label{eq:kdpp_map}
	Y^\star = \argmax_{Y \subseteq \mc{Y},\, |Y| = r} \mc{P}_{\m{L}}^r(Y),
\end{equation}
where $\mc{P}_{\m{L}}^r(Y) = \frac{\det(\m{L}_Y)}{e_r(\lambda_1, \ldots, \lambda_n)}$ and $e_r(\lambda_1, \ldots, \lambda_n) = \sum_{|Y'| = r} \det(\m{L}_{Y'})$ is the $r$-th elementary symmetric polynomial of the eigenvalues of $\m{L}$. Then solving~\eqref{eq:kdpp_map} is equivalent to solving \probref{prob:node_selection}.
\end{myprs}

\begin{proof}
For fixed $r$, the denominator $e_r(\lambda_1,\ldots,\lambda_n)$ is a positive constant independent of $Y$; therefore maximizing $\mc{P}_{\m{L}}^r(Y)$ is equivalent to maximizing $\det(\m{L}_Y)$. Since $\log(\cdot)$ is strictly increasing on $\mbb{R}_{>0}$ and $\m{L}=\m{W}_o\succ\m{0}$, this is equivalent to maximizing $\logdet{(\m{W}_o)_Y}$ in~\eqref{eq:sns_problem}. This completes the proof.
\end{proof}

The $k$-DPP is the cardinality constrained counterpart of the L-ensemble (\defref{def:L_ensemble}); it preserves the repulsive structure of the DPP while restricting to subsets of exactly $r$ nodes. By~\propref{pro:kdpp_map}, maximizing $\mc{P}_{\m{L}}^r(Y)$ over subsets of size $r$ reduces to maximizing $\logdet{(\m{W}_o)_Y}$. With $\m{C}=\m{I}$, the Gramian satisfies $\m{L}=\m{W}_o\succeq\m{I}$. Since $\logdet{\m{L}_Y}$ is submodular and monotone nondecreasing with $\logdet{\m{L}_\emptyset} = 0$ (\lemref{lem:submod}), the greedy algorithm that sequentially adds the node with the largest marginal gain yields a $(1-1/e)$ approximation~\cite{Nemhauser1978}. The greedy selection is summarized in \algref{alg:greedy_map} and returns a single deterministic configuration that approximates the MAP of the $k$-DPP.
At each iteration, the marginal gain of a candidate node is evaluated through the Schur complement $\m{L}_{jj} - \m{L}_{j,Y}\m{L}_Y^{-1}\m{L}_{Y,j}$ of the Gramian, which is the residual observability energy of node $j$ that is orthogonal to the observability directions of the selected set; the gain of every candidate decreases as correlated nodes are added. This is consistent with the conditional inclusion probability of \coref{cor:cond_incl}.

\begin{algorithm}[t]
	\caption{$k$-DPP sampling for node selection.}\label{alg:kdpp}
	\DontPrintSemicolon
	\textbf{Input:} DPP kernel $\m{L} = \m{W}_o \in \Rn{n \times n}$, subset cardinality $r$\;\vspace{0.1cm}
	\textbf{Compute:} eigendecomposition $\m{L} = \m{V}\m{\Lambda}\m{V}^\top$, $\m{\Lambda} = \mr{diag}(\lambda_1,\ldots,\lambda_n)$\;\vspace{0.1cm}
	// \textit{select $r$ eigenvectors via elementary symmetric polynomial}\;
	Compute $e_l^t$ for $l = 0,\ldots,r$, $t = 1,\ldots,n$ via the recursion $e_l^t = e_l^{t-1} + \lambda_t\, e_{l-1}^{t-1}$\;
	$J \leftarrow \emptyset$,\; $l \leftarrow r$\;
	\For{$t \leftarrow n$ \KwTo $1$}{
		\If{$u \sim \mr{Uniform}(0,1) < \lambda_t\, e_{l-1}^{t-1} \,/\, e_l^t$}{
			$J \leftarrow J \cup \{t\}$,\quad $l \leftarrow l - 1$\;
		}
		\lIf{$l = 0$}{\textbf{break}}
	}
	// \textit{sample nodes from selected eigenvector subspace}\;\vspace{0.1cm}
	$\m{V}_J \leftarrow [\m{v}_j]_{j \in J} \in \Rn{n \times r}$\;
	$Y \leftarrow \emptyset$\;
	\For{$s \leftarrow r$ \KwTo $1$}{
		\textbf{select:} node $i$ with probability $\prob(i) = \tfrac{1}{s} \|\m{V}_J^\top \m{e}_i\|^2$\;
		$Y \leftarrow Y \cup \{i\}$\;
		\textbf{update:} $\m{V}_J \leftarrow \mr{orth}\!\left(\m{V}_J - \tfrac{\m{V}_J\m{V}_J^\top \m{e}_i\, \m{e}_i^\top \m{V}_J}{\|\m{V}_J^\top \m{e}_i\|^2}\right)$
		// \textit{project out $\m{e}_i$ and orthogonalize the basis}\;
	}
	\textbf{Output:} node configuration $Y$ with $|Y| = r$\;
\end{algorithm}

While a single MAP solution can be obtained greedily, the same formulation yields a distribution over all subsets of size $r$ through the $k$-DPP. This enables sampling many feasible sensor configurations, where each sample is generally suboptimal relative to the MAP set, yet is weighted by $\det(\m{L}_Y)$ and therefore concentrates on diverse subsets. Exact $k$-DPP sampling runs in $\mc{O}(n^3 + nr^2)$ time via the spectral algorithm of~\cite{Hough2006,Kulesza2012}; it computes the eigendecomposition of $\m{L}$, selects $r$ observability modes through the elementary symmetric polynomial recursion, and sequentially samples nodes while projecting out the directions already selected. \algref{alg:kdpp} is based on the spectral sampling algorithm of~\cite{Hough2006,Kulesza2012} and applied to the Gramian $\m{L} = \m{W}_o$, where every sample is a node configuration weighted by its observability volume (\theoref{thm:gramian_dpp}). The resulting ensemble is used in \secref{sec:case_study} for state estimation under node failures, where $k$-DPP samples replace nodes that fail, without solving \probref{prob:node_selection} another time. This, along with the greedy guarantee for a single MAP solution (\lemref{lem:submod}), provides a computationally tractable way to generate an ensemble of diverse candidate sensor configurations. When the marginal kernel in~\eqref{eq:marginal_K} is parameterized through a structured low rank approximation, DPP inference and parameter learning can be performed with complexity that can be sublinear in the number of nodes~\cite{Dupuy2018}. For large scale systems with $\m{K}=\m{W}_o(\m{W}_o+\m{I})^{-1}$, a small effective observable rank in~\eqref{eq:expected_card} indicates that a few observability modes dominate, which supports low rank approximations of $\m{W}_o$ for computations that are sublinear in $n$ and efficient solution of~\eqref{eq:kdpp_map} with diverse sensor configurations.

Equation~\eqref{eq:kdpp_map} also admits a continuous relaxation through the multilinear extension $F(\m{x}) := \log \det(\mr{diag}(\m{x})(\m{L} - \m{I}) + \m{I})$, where $\m{x} \in [0,1]^n$ and $F(\m{1}_Y) = \log\det(\m{L}_Y)$ at every vertex~\cite{Gillenwater2012}. Under the cardinality constraint $\sum_i x_i = r$, Frank-Wolfe optimization converges to a stationary point, where rounding yields a $1/4$-approximation for general DPPs~\cite{Gillenwater2012}. The multilinear relaxation can be used to solve the MAP problem in~\eqref{eq:kdpp_map} as shown in~\cite{Kazma2026a}; it is useful when constraints beyond cardinality are imposed.

\section{Case Study}\label{sec:case_study}
To demonstrate and validate the theoretical results established in the aforementioned sections, we investigate the following questions.
\begin{itemize}[leftmargin=*]
		\item$(\mr{Q}1)$ Does $k$-DPP sampling~\eqref{eq:kdpp_map} via \algref{alg:kdpp} produce diverse node configurations efficiently while retaining large $\logdett$ values compared to MAP greedy optimization of~\eqref{eq:sns_problem} via \algref{alg:greedy_map}?
		\item$(\mr{Q}2)$ Do the marginal kernel (\theoref{thm:marginal_modal}), the inclusion monotonicity (\coref{cor:stoch_dom}), and the quality and diversity decomposition~\eqref{eq:quality_diversity} with its negative dependence (\coref{cor:neg_dep}) capture the observability structure of the networks?
		\item$(\mr{Q}3)$ Does the $\logdett$ ordering of node configurations carry over to the estimation error of \propref{prop:estimation}, and do the $k$-DPP samples attain lower errors than configurations drawn uniformly at random? When selected nodes fail, do the $k$-DPP samples yield configurations of $r$ nodes that avoid the failed ones without solving \probref{prob:node_selection} another time?
		\item$(\mr{Q}4)$ Do the results extend to linear output matrices (\theoref{thm:general_C}) and to actuator node selection through the controllability Gramian of \secref{subsec:control}?
\end{itemize}

Questions $(\mr{Q}1)$ and $(\mr{Q}3)$ evaluate the node selection performance of the framework, whereas $(\mr{Q}2)$ and $(\mr{Q}4)$ validate the derived properties; each question is addressed in one of the following subsections. We consider three LTI network models of increasing complexity studied in the literature~\cite{Liu2011}. The first is a directed small world network, denoted as directed Watts-Strogatz (DWS), constructed from the undirected model of~\cite{Watts1998} with $n = 122$ nodes, in which each node is coupled to its four nearest neighbors, with rewiring probability $p_{\mr{rew}} = 0.15$. The undirected graph is converted to a digraph and $40\%$ of directed edges are randomly pruned to produce a sparse directed network with small world connectivity. The second and third are canonical random graph models, \textit{(i)} an Erd\H{o}s-R\'{e}nyi (ER) random graph~\cite{Erdos1959} with $n = 200$ nodes and edge probability $p_{\mr{edge}} = 0.05$; and \textit{(ii)} a Barab\'{a}si-Albert (BA) preferential attachment graph~\cite{Barabasi1999} with $n = 200$ nodes and $m_0 = 3$ edges per new node. Both of the random graph models are used in~\cite{Liu2011} to study the controllability of networks. 

We generate $20$ independent random realizations for each of the three network models. For all three networks, we consider the canonical linear dynamics $\m{x}[k+1] = \m{A}\m{x}[k]$~\cite{Pasqualetti2014a}, where $\m{A}$ is obtained from the adjacency matrix of the network and scaled as $\m{A} \leftarrow \alpha \m{A} / \rho(\m{A})$ with $\alpha = 0.95$, which gives $\rho(\m{A}) = \alpha < 1$ and therefore ensures Schur stability~\cite{Pasqualetti2014a}. We set $\m{C} = \m{I}_n$, where each node is a candidate sensor. The observation horizon is $N = 200$, and the sensor budget is $r \in \{5, 10, 15, 20\}$. For each network and each $r$, we draw $M = 500$ independent samples from the $k$-DPP using \algref{alg:kdpp}. We compare these samples against two deterministic node selection methods, \textit{(i)} greedy under $\logdet{(\m{W}_o(Y) + 10^{-8}\,\m{I})}$, where the small regularization follows~\cite{Summers2016} and allows for $\m{W}_o(Y)$ to be full rank; and \textit{(ii)}  MAP greedy (\algref{alg:greedy_map}) evaluated under the objective of \probref{prob:node_selection}. For the marginal kernel validation in $(\mr{Q}2)$, we also draw $M_{\mr{val}} = 200$ samples from the unconstrained DPP for each network. %

\begin{figure}[t]
	\centering
	\includegraphics[width=\linewidth]{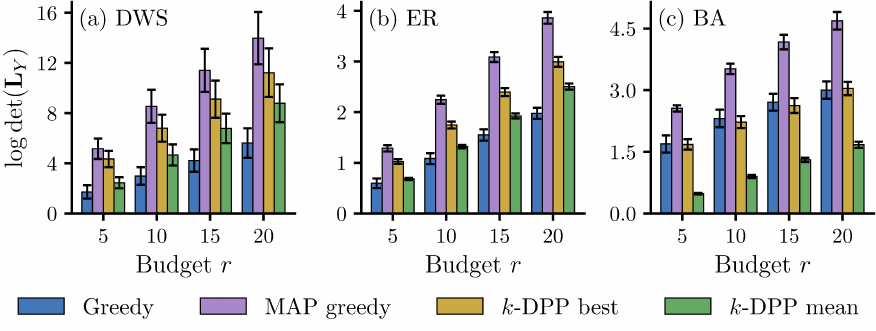}
\caption{$k$-DPP sampling and MAP greedy for $\logdet{\m{L}_Y}$ across sensor budgets $r \in \{5, 10, 15, 20\}$. For each network and budget, we depict greedy under $\logdet{\m{W}_o(Y)}$, MAP greedy, the best $k$-DPP sample, and the mean $k$-DPP sample over $M = 500$ draws. Bars show means over $20$ system realizations; error bars indicate $\pm 1$ standard deviation across realizations. The $k$-DPP samples achieve $\logdett$ ratios $\varrho \in [0.63, 0.84]$ relative to the MAP greedy value.} \label{fig:logdet_comparison}
\end{figure}

\subsection{Node selection under MAP}\label{subsec:cs_performance}
\figref{fig:logdet_comparison} depicts the DPP objective $\logdet{\m{L}_Y}$ for MAP greedy (\algref{alg:greedy_map}), the greedy selection under $\logdet{\m{W}_o(Y)}$, and $k$-DPP sampling (\algref{alg:kdpp}) across all networks and budgets $r$, including the best and mean values over $M = 500$ independent $k$-DPP samples. The greedy selection under $\logdet{\m{W}_o(Y)}$ results in a lower objective value and lies below the mean of the $k$-DPP samples for the DWS and ER networks. This gap follows from the rank deficiency of the parameterized Gramian at the considered budgets. For $r = 10$, the selected configurations attain a rank of $\m{W}_o(Y)$ of at most $101$ of the $122$ states for the DWS network and at most $77$ of the $200$ states for the ER network across all realizations, and thus the regularized objective is dominated by the unobserved directions. In contrast, every principal submatrix of $\m{W}_o \succ \m{0}$ is positive definite, and thus \probref{prob:node_selection} remains well-posed at every budget without regularization. The DPP formulation therefore allows considering budgets where $\logdet{\m{W}_o(Y)}$ is degenerate. We note that single swap local search~\cite{Gillenwater2012} does not improve upon \algref{alg:greedy_map} in any instance. 

\begin{figure}[t]
	\centering
	\includegraphics[width=\linewidth]{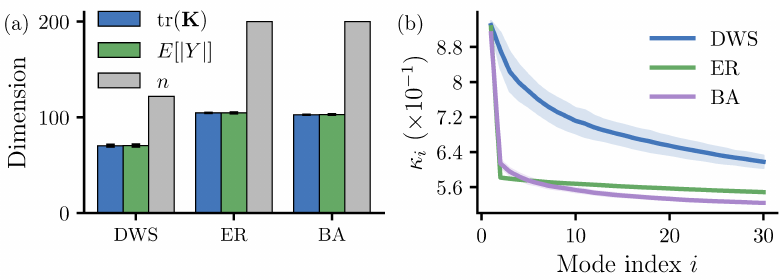}
	\vspace{-0.3cm}
	\caption{(a)~Effective observable rank validation, theoretical value $n_{\mr{eff}} = \trace(\m{K})$ is equivalent to the empirical mean cardinality $\mbb{E}[|Y|]$ over $M_{\mr{val}} = 200$ DPP samples. (b)~Observability mode activation probabilities $\kappa_i = \lambda_i / (\lambda_i + 1)$. Bars and curves show means over $20$ realizations; bands indicate $\pm 1$ standard deviation.}
	\label{fig:marginal_kernel}
\end{figure}

The $k$-DPP samples achieve $\logdett$ ratios $\varrho \in [0.63, 0.84]$ relative to the MAP greedy value, with all $500$ samples being unique configurations that collectively include $100\%$ of the network nodes and have mean pairwise Jaccard distance $\bar{d}_J > 0.90$, where $\bar{d}_J(Y_i, Y_j) = 1 - |Y_i \cap Y_j|/|Y_i \cup Y_j|$ measures the fraction of nodes that differ between two samples. For $M = 500$, the DPP therefore provides distinct configurations with large $\logdett$ values while including every node across the samples. The DWS network exhibits $\logdet{\m{L}_Y}$ standard deviations that are at least $10\times$ larger than ER and BA ($1.31$ vs. $0.08$ and $0.13$ at $r = 10$). This follows from the asymmetric edge pruning step that amplifies topological variability across realizations. The $\logdett$ ratios $\varrho$, however, remain consistent, with a standard deviation of at most $0.06$. This confirms that the DPP approximation is robust to the network realization. The per sample $k$-DPP time is $0.001$-$0.010$~s, with the $\mc{O}(n^3)$ eigendecomposition computed once per network; it is comparable to the greedy time of $0.004$-$0.024$~s. The above answers $(\mr{Q}1)$ and shows how the proposed framework enables a probabilistic approach to node selection.

\subsection{DPP diversity and marginal kernel results}\label{subsec:cs_diversity}
The joint inclusion probability $\prob(i, j \in Y) = \m{K}_{ii}\m{K}_{jj} - \m{K}_{ij}^2$ (\coref{cor:neg_dep}) is validated under standard DPP sampling with $2000$ samples per network, where the root mean square error (RMSE) between the empirical and theoretical values is at most $0.012$. This confirms that the DPP suppresses redundant node selection by enforcing negative dependence, as the joint inclusion probability is always less than or equal to the product of the marginal inclusion probabilities $\prob(i \in Y) \cdot \prob(j \in Y)$. The marginal kernel $\m{K}$ reveals the observability structure through \figref{fig:marginal_kernel}. \textit{(i)}~The cardinality identity $\mbb{E}[|Y|] = \trace(\m{K})$ (\theoref{thm:marginal_modal}) is validated empirically, with theoretical values of $70.3$, $104.7$, and $102.7$ for DWS, ER, and BA and empirical means $70.4$, $104.7$, and $102.9$ over $M_{\mr{val}}=200$ samples. Since the network dimensions differ, the corresponding normalized values $n_{\mr{eff}}/n$ are $0.576$, $0.524$, and $0.514$, respectively.
 \textit{(ii)}~The modal probabilities $\kappa_i$ (\theoref{thm:marginal_modal}) decay gradually for the ER network while showing a rapid initial drop followed by a long tail for the BA network. Since each $\kappa_i = \lambda_i/(\lambda_i + 1)$ is the inclusion probability of the $i$-th observability mode, this indicates that the BA network concentrates most of its Gramian energy in a few dominant modes while many modes contribute $\kappa_i \approx 0$, consistent with its community structure topology. These results address the marginal kernel and negative dependence questions posed in $(\mr{Q}2)$.

\begin{figure}[t]
	\centering
	\includegraphics[width=\linewidth]{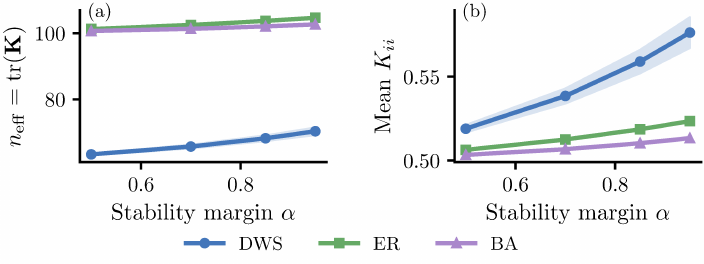}
	\vspace{-0.3cm}
	\caption{Inclusion monotonicity validation (\coref{cor:stoch_dom}). As the stability margin $\alpha$ increases from $0.5$ to $0.95$, the Gramian $\m{W}_o$ increases. This results in monotone increases in \textit{(a)} the effective observable rank $n_{\mr{eff}} = \trace(\m{K})$ and \textit{(b)} the mean marginal inclusion probability $\bar{\m{K}}_{ii}$. Lines show means over $20$ realizations; bands indicate $\pm 1$ standard deviation.}
	\label{fig:inclusion_monotonicity}
\end{figure}

With that in mind, we now vary the stability margin $\alpha \in \{0.5, 0.7, 0.85, 0.95\}$, thereby changing $\m{A} \leftarrow \alpha \m{A} / \rho(\m{A})$, such that $\rho(\m{A}) = \alpha < 1$. For values of $\alpha_1 < \alpha_2$, this implies that $\m{W}_o^{(\alpha_1)} \preceq \m{W}_o^{(\alpha_2)}$. According to \coref{cor:stoch_dom}, we obtain $\m{K}^{(\alpha_1)} \preceq \m{K}^{(\alpha_2)}$; this is depicted in \figref{fig:inclusion_monotonicity}. The effective rank $n_{\mr{eff}}$ increases monotonically with $\alpha$ for all networks (from $63.3$ to $70.3$ for DWS, $101.3$ to $104.7$ for ER, $100.7$ to $102.7$ for BA), with the DWS network exhibiting the steepest increase. Since larger $\alpha$ amplifies the Gramian eigenvalues $\lambda_i$, the mapping $\kappa_i = \lambda_i/(\lambda_i + 1)$ pushes more modes toward full activation ($\kappa_i \approx 1$), increasing $n_{\mr{eff}}$ and raising the inclusion probability of every sensor configuration; see \coref{cor:stoch_dom}. This behavior follows from the slower decay of $\m{A}^k$ at a larger stability margin. The outputs accumulate more observability energy over the horizon, and additional modes count toward the effective observable rank of \secref{subsec:marginal}. The DWS network is most sensitive to this effect because its smaller dimension ($n = 122$) and sparse directed topology concentrate Gramian energy in fewer modes, so the transition from weakly to strongly activated modes is steeper. This empirically validates the inclusion monotonicity part of $(\mr{Q}2)$.

\begin{figure}[t]
	\centering
	\includegraphics[width=\linewidth]{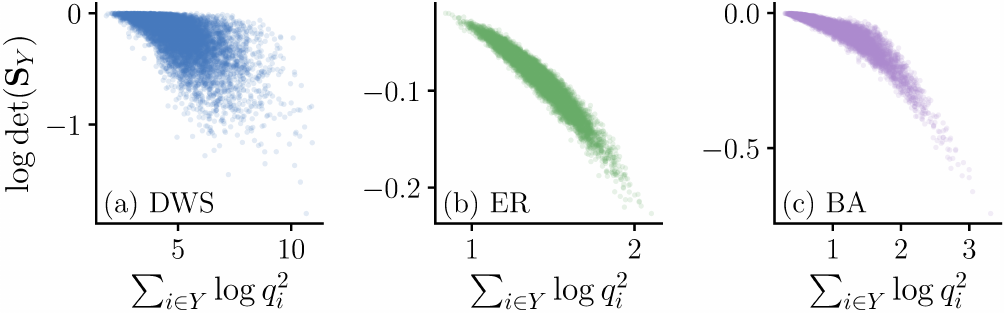}
	\vspace{-0.3cm}
	\caption{For $M = 500$ $k$-DPP samples at $r = 10$ and for all $20$ realizations, each point represents one sample. The $x$-axis is the quality term $\sum_{i \in Y} \log q_i^2$ and the $y$-axis is the diversity term $\logdet{\m{S}_Y}$.}
	\label{fig:quality_diversity}
\end{figure}

\figref{fig:quality_diversity} depicts the decomposition~\eqref{eq:quality_diversity} for $M = 500$ $k$-DPP samples for $r = 10$. A negative correlation between quality $\sum_{i \in Y} \log q_i^2$, where $q_i^2 = (\m{W}_o)_{ii}$ is the individual observability energy of node $i$, and diversity $\logdet{\m{S}_Y}$, where $\det(\m{S}_Y)$ measures the orthogonality of the selected observability directions, is evident across all networks. Configurations that select nodes with high individual Gramian energy tend to have more collinear observability vectors, reducing $\det(\m{S}_Y)$. This tradeoff is pronounced in the BA network, where nodes forming communities share similar Gramian rows. The ER and BA networks exhibit correlations of $-0.98$ and $-0.92$, whereas the DWS network yields a weaker correlation of $-0.69$, consistent with its sparse directed topology. That is, the DWS network yields a Gramian with a structure that varies more across nodes. This confirms that the $\logdett$ metric implicitly balances these two objectives, as established in~\eqref{eq:quality_diversity}.

\subsection{State estimation under node selection and failures}\label{subsec:cs_estimation}
We consider more heterogeneous dynamics over the three networks. The diagonal entries of $\m{A}$ and the coupling weights on the directed edges are drawn uniformly and independently from $[0.05, 0.95]$ and $[0.1, 1]$ for each realization, and the coupling term is scaled such that $\rho(\m{A}) = 0.99$. The heterogeneous dynamics and the larger stability margin increase the Gramian (\coref{cor:stoch_dom}) and differentiate the observability energies $(\m{W}_o)_{ii}$ of the nodes, and thus the node configurations are distinguishable through their estimation errors. By \propref{prop:estimation}, the conditional covariance of the estimation errors associated with the selected state variables is $\sigma^2 ((\m{W}_o)_Y)^{-1}$ for measurements collected over the horizon with noise of variance $\sigma^2$. Maximizing the objective of \probref{prob:node_selection} minimizes the determinant of this covariance. We therefore compute the estimation error as
\begin{equation}\label{eq:cond_error}
	e(Y) := \sqrt{\trace\big(\sigma^2 ((\m{W}_o)_Y)^{-1}\big) / |Y|},
\end{equation}
which is obtained from the covariance~\eqref{eq:cond_cov} and is the root mean square of the estimation error of the selected state variables conditioned on the remaining ones and normalized by $|Y|$, and we set $\sigma^2 = 1$. The error~\eqref{eq:cond_error} is computed in closed form. This is not the error of an estimate of the full state from the sensors in $Y$; such an estimate might not be defined at the considered budgets. This is due to $\m{W}_o(Y)$ being rank deficient (\secref{subsec:cs_performance}), where the error is unbounded along the unobserved directions. Equation~\eqref{eq:cond_error} instead quantifies the accuracy with which the selected state variables are estimated under the measurements of \propref{prop:estimation}, and it is well-defined for every budget. We note that the squared error in~\eqref{eq:cond_error} is the normalized trace of the covariance in~\eqref{eq:cond_cov}, whereas the objective of \probref{prob:node_selection} is the determinant of its inverse. The two are thus distinct measures of the same covariance.

\begin{figure}[t]
	\centering
	\includegraphics[width=\linewidth]{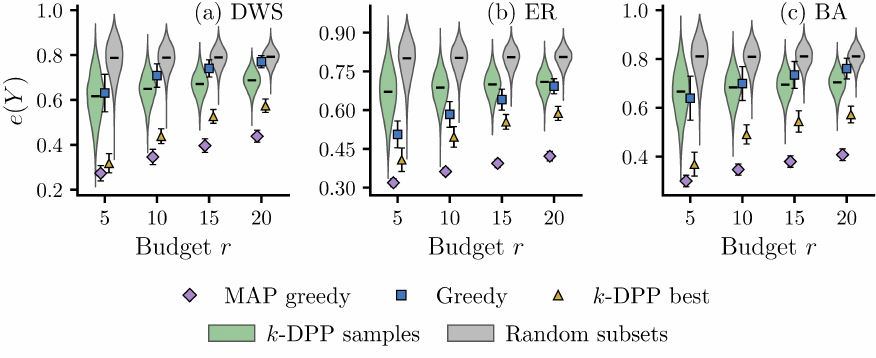}
	\vspace{-0.3cm}
	\caption{State estimation error~\eqref{eq:cond_error} for the three networks. The shaded regions depict the distribution of the error over the $M = 500$ $k$-DPP samples (\algref{alg:kdpp}) and the $M = 500$ configurations drawn uniformly at random from all $20$ realizations; horizontal lines indicate medians. Markers depict the MAP greedy (\algref{alg:greedy_map}), greedy, and the largest $\logdett$ $k$-DPP sample; error bars indicate $\pm 1$ standard deviation across realizations.}
	\label{fig:estimation_comparison}
	\vspace{0.1cm}
\end{figure}

\figref{fig:estimation_comparison} depicts the distribution of the error~\eqref{eq:cond_error} over the $M = 500$ $k$-DPP samples and over $M = 500$ configurations drawn uniformly at random from all $20$ realizations. The figure also shows the errors of the two greedy configurations and of the highest $\logdett$ $k$-DPP sample. For every network and every budget, configurations with larger $\logdett$ values attain lower estimation errors. At $r = 10$, the errors for the DWS network are $0.346$ for the MAP greedy configuration, $0.439$ for the highest $\logdett$ sample, and $0.648$ and $0.786$ for the mean over the $k$-DPP samples and the uniformly drawn configurations. A $k$-DPP sample results in a lower error than a randomly drawn configuration of the same cardinality with probability $0.85$-$0.90$ at $r = 5$ and $0.97$-$0.98$ at $r = 20$ across the three networks. The greedy selection under $\logdet{\m{W}_o(Y)}$ results in estimation errors of $0.708$, $0.584$, and $0.699$ at $r = 10$ for the DWS, ER, and BA networks, which are similar to the mean over the $k$-DPP samples and lower than the uniformly drawn configurations. The $\logdett$ value of a configuration is therefore indicative of its estimation error; this answers the estimation part of $(\mr{Q}3)$. Furthermore, we note that the error~\eqref{eq:cond_error} increases with the budget $r$ for every method, whereas the error of an estimate decreases when additional measurements are collected. Increasing $r$ adds a state variable to the estimated set in~\eqref{eq:cond_cov}, where its estimate is affected by the measurement noise. The added variable carries a smaller residual observability energy than the nodes selected before it (\coref{cor:cond_incl}), and thus the normalized error per selected variable increases. 

\begin{figure}[t]
	\centering
	\includegraphics[width=\linewidth]{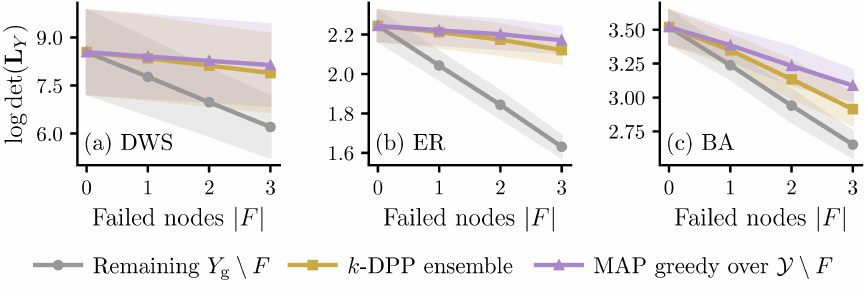}
	\vspace{-0.3cm}
	\caption{Node failure experiment at $r = 10$, where the nodes $F$ of the MAP greedy configuration $Y_{\mr{g}}$ fail, $|F| \in \{1,2,3\}$ ($10$ draws of $F$ per value of $|F|$). The three strategies measure the remaining nodes $Y_{\mr{g}} \setminus F$, restore the budget through $Y_{\mr{res}}$ with nodes from the stored samples ($k$-DPP ensemble), and apply \algref{alg:greedy_map} over the remaining ground set $\mc{Y} \setminus F$. Lines show means over $20$ realizations; bands indicate $\pm 1$ standard deviation.}
	\label{fig:sensor_failure}
\end{figure}

We now consider the setting in which $r$ nodes are active and a subset of these nodes remains active at a time due to sensing budgets or sensor failures. Changing the active configuration therefore requires selecting a different node subset. Let $Y_{\mr{g}} \subseteq \mc{Y}$ denote the MAP greedy configuration (\algref{alg:greedy_map}) at $r = 10$, and let $F \subset Y_{\mr{g}}$ denote the set of failed nodes, where a node in $F$ can no longer be measured. We consider $|F| \in \{1, 2, 3\}$ over $10$ random draws of $F$ per realization, and we compare three strategies. The first measures the remaining nodes $Y_{\mr{g}} \setminus F$ of the greedy configuration, with cardinality $r - |F|$. The second retains the nodes $Y_{\mr{g}} \setminus F$ and identifies the $k$-DPP sample $Y_{\mr{s}}$ with the largest $\logdett$ value among the stored samples that satisfy $Y_m \cap F = \emptyset$. This strategy restores the budget through the configuration $Y_{\mr{res}} := (Y_{\mr{g}} \setminus F) \cup A$, where $A \subseteq Y_{\mr{s}} \setminus Y_{\mr{g}}$ collects the $|F|$ nodes with the largest marginal gains with respect to the retained nodes according to~\coref{cor:cond_incl}, in which the conditioning enters through the Schur complement of the kernel. Given that $Y_{\mr{g}} \setminus F \subset Y_{\mr{res}}$ and $\logdet{\m{L}_Y}$ is monotone (\lemref{lem:submod}), the second strategy improves on the first in every draw. The third applies \algref{alg:greedy_map} over the remaining ground set $\mc{Y} \setminus F$, which returns a new configuration of $r$ nodes.~\figref{fig:sensor_failure} depicts the resulting $\logdet{\m{L}_Y}$. Measuring $Y_{\mr{g}} \setminus F$ loses observability as $|F|$ increases (from $8.54$ to $6.20$ for DWS). Restoring the budget with $Y_{\mr{res}}$ results in values $8.35$, $8.12$, and $7.89$ for $|F| = 1, 2, 3$ on the DWS network. Applying MAP greedy over $\mc{Y} \setminus F$ yields the highest $\logdett$ value ($8.40$, $8.27$, and $8.14$), at the cost of solving \probref{prob:node_selection} after every failure event. That being said, the restored configurations recover at least $94\%$ of this value across the three networks and failures without resolving \probref{prob:node_selection}, where $66\%$-$95\%$ of the $M = 500$ samples satisfy $Y_m \cap F = \emptyset$ and the added nodes are thus selected from samples that are already available. The stored $k$-DPP samples therefore restore configurations of $r$ nodes that avoid the failed ones without solving \probref{prob:node_selection}. This offers a probabilistic alternative to solving the sensor selection problem, thereby answering $(\mr{Q}3)$.
\begin{figure}[t]
	\centering
	\includegraphics[width=\linewidth]{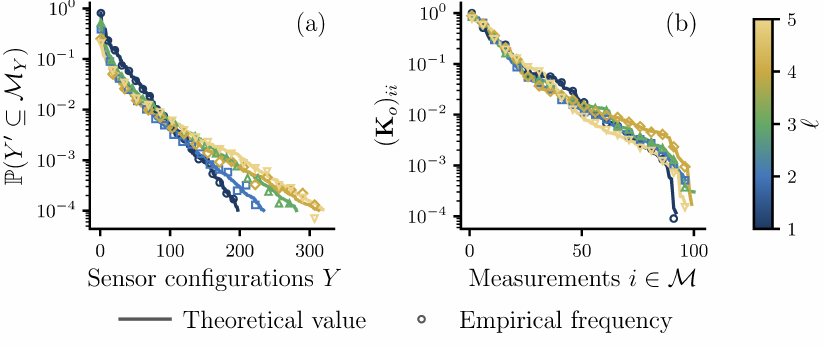}
	\vspace{-0.3cm}
	\caption{Validation of \theoref{thm:general_C} for five output matrices. Curves depict the theoretical values in decreasing order, and markers depict the corresponding empirical frequencies over $M = 10^5$ samples of the $k$-DPP (\algref{alg:kdpp}). (a)~The probability $\det(\m{W}_o(Y))/\det(\m{W}_o)$ of~\eqref{eq:sensor_incl} over all sensor configurations of cardinality $|Y| \in \{6,7,8\}$. (b)~The marginal inclusion probability $(\m{K}_o)_{ii}$ over all measurements $i \in \mc{M}$.}
	\label{fig:general_C_validation}
\end{figure}
\vspace{-0.4cm}
\subsection{Linear output matrices and actuator node selection}\label{subsec:cs_general}
We consider a system with $n = 10$ states obtained from the DWS model and $p = 10$ candidate sensors. We construct five output matrices $\m{C}$ in which each sensor measures $\ell \in \{1,\ldots,5\}$ state variables with random weights. For each output matrix we draw $M = 10^5$ samples of the $k$-DPP of the kernel $\m{L}_{\mc{M}}$ using \algref{alg:kdpp}, and consider all sensor configurations of cardinality $|Y| \in \{6,7,8\}$. \figref{fig:general_C_validation} depicts the theoretical values sorted in decreasing order. The corresponding empirical frequencies are depicted as markers. The empirical frequencies follow the theoretical curves over more than three orders of magnitude in probability. Thus the agreement between the theoretical and empirical values holds for considered output mappings.

That being said, \theoref{thm:general_C} holds for all five output matrices. The nonzero eigenvalues of $\m{L}_{\mc{M}}$ match the eigenvalues of $\m{W}_o$ to within $1.9 \times 10^{-13}$, and the normalization satisfies $|\log e_n - \logdet{\m{W}_o}| \leq 5.0 \times 10^{-14}$. The kernel~\eqref{eq:proj_kernel} satisfies $\trace(\m{K}_o) = n$, and the projection identity $\m{K}_o^2 = \m{K}_o$ holds with an error of at most $8.3 \times 10^{-15}$. Note that the event probabilities depend on the output matrix. For $\ell = 1$ the probability is concentrated on a few configurations, whereas denser output matrices spread the probability across more configurations. Furthermore, for each output matrix, the empirical frequency with which the sampled measurements are all generated by the sensors in $Y$ agrees with $\det(\m{W}_o(Y))/\det(\m{W}_o)$ with RMSE at most $3.8 \times 10^{-4}$ over all $375$ configurations of the considered cardinalities. The marginal inclusion frequency of each measurement matches $(\m{K}_o)_{ii}$ with RMSE at most $8.0 \times 10^{-4}$. For every output matrix and cardinality, the configuration that maximizes $\logdet{\m{W}_o(Y)}$ has the largest probability of the event $Y' \subseteq \mc{M}_Y$ among configurations of the same cardinality.

\begin{figure}[t]
	\centering
	\includegraphics[width=\linewidth]{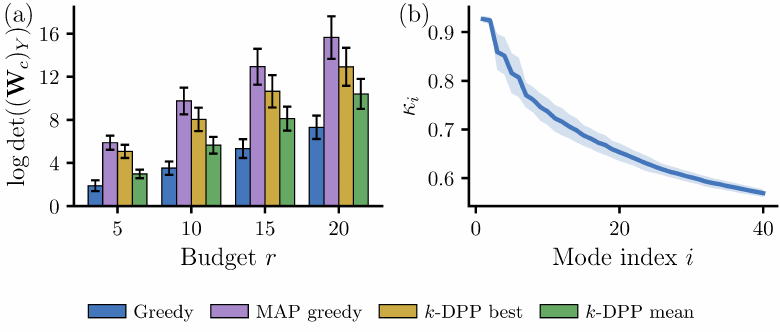}
	\vspace{-0.3cm}
	\caption{Actuator node selection on the directed lattice network using the controllability Gramian $\m{W}_c$. (a)~$\logdet{(\m{W}_c)_Y}$ for MAP greedy, and the best and mean $k$-DPP samples over $M = 500$ draws. (b)~Mode activation probabilities $\kappa_i$ of the controllability Gramian. Error bars and bands indicate $\pm 1$ standard deviation across $20$ realizations.}
	\label{fig:actuator}
\end{figure}

We now consider a directed lattice network that represents a spatially distributed process, constructed from a $10 \times 10$ grid graph converted to a digraph with $30\%$ of the directed edges randomly pruned. Actuator placement on network models representing a grid graph is studied in~\cite{Summers2016}. Each node is a candidate actuator, such that $\m{B} = \m{I}_n$, and by the duality of \secref{subsec:control} the DPP matrix is the controllability Gramian $\m{L} = \m{W}_c = \m{W}_o(\m{A}^\top, \m{I})$. We generate $20$ independent realizations of the edge pruning, and the horizon is $N = 200$. We solve the actuator node selection problem by maximizing $\logdet{(\m{W}_c)_Y}$ over configurations of cardinality $r \in \{5,10,15,20\}$ using MAP greedy (\algref{alg:greedy_map}) and greedy, and we draw $M = 500$ samples using \algref{alg:kdpp}. \figref{fig:actuator} depicts the results. The $k$-DPP samples attain $\logdett$ ratios $\varrho \in [0.77, 0.95]$ relative to the greedy value across budgets and realizations, which is consistent with the ratios obtained for sensor node selection in \secref{subsec:cs_performance}. The mode activation probabilities $\kappa_i$ of the controllability Gramian decay gradually, and the effective rank is $n_{\mr{eff}} = 58.4$ on average. The framework therefore extends to actuator node selection. The above results answer $(\mr{Q}4)$ and conclude the case studies.

\section{Paper Summary, Limitations and Future Work}\label{sec:conclusion}
This paper establishes connections between Gramians in control and DPPs. When each state variable is a candidate node, the observability Gramian of a discrete-time LTI system is a valid DPP L-ensemble matrix over the node set, and the principal submatrix node selection problem admits a $k$-DPP MAP reformulation. For a linear output matrix, the determinant of the Gramian associated with a sensor subset is proportional to the probability that a DPP selects only rows of the observability matrix associated with those sensors. Based on these connections, we derive a marginal kernel result, an effective observable rank that quantifies how many observability modes are activated for a given system, a decomposition of the DPP probability that separates individual node observability energy from node diversity, and an inclusion monotonicity result. We validate the results on several networks and show that the DPP samples provide configurations that can replace failed sensor nodes without solving the node selection problem.

This work is not devoid of limitations. Specifically, we do not consider: \textit{(i)} process or measurement noise in the problem formulation, although we study state estimation performance under the obtained configurations; \textit{(ii)} node selection under constraints beyond a cardinality budget, such as partition matroid or communication constraints; and \textit{(iii)} models that depict infrastructure system networks, as the case studies consider network topologies studied in the literature. The work presented thus merits future investigation on this topic, which entails: \textit{(a)} addressing the aforementioned limitations; \textit{(b)} investigating whether empirical and variational Gramians for time-varying and nonlinear systems yield results analogous to those established here; and \textit{(c)} learning the DPP matrix from measurement data rather than deriving it from a known model.

\appendices
\section{DPPs for Linear Output Matrices}\label{app:generalC}
For a linear output matrix $\m{C} \in \Rn{p \times n}$, each sensor node measures a linear combination of one or multiple state variables, rather than a single state variable. Node selection over the $p$ candidate sensors is parameterized by $\bm{\gamma} \in \{0,1\}^p$ through $\m{W}_o(\bm{\gamma}) = \sum_{i=1}^{p} \gamma_i \m{W}_o^{(i)}$, where $\m{W}_o^{(i)} = \sum_{k=0}^{N-1} (\m{A}^k)^\top \m{c}_i \m{c}_i^\top \m{A}^k$ and $\m{W}_o(\bm{\gamma}) = \m{W}_o(Y)$ for $Y = \{i \mid \gamma_i = 1\}$. By the Cauchy-Binet formula, $\det(\m{W}_o(Y))$ is a sum of squared minors of $\m{O}_Y$ rather than a principal minor of $\m{W}_o$, thus the L-ensemble $\m{L} = \m{W}_o$ does not directly represent this sensor selection. This selection is characterized by \theoref{thm:general_C} through a DPP over the rows of the observability matrix.

Let the ground set $\mc{M} = \{1, \ldots, Np\}$ index the rows of $\m{O}$, where $\tilde{\m{o}}_i^\top$ denotes the $i$-th row and $\m{O}_S$ collects the rows indexed by $S \subseteq \mc{M}$. For a sensor configuration $Y \subseteq \mc{Y}$, the set $\mc{M}_Y \subseteq \mc{M}$ collects the rows generated by the sensors in $Y$ over the horizon $N$; these rows form $\m{O}_Y$ in~\eqref{eq:obs_matrix_row}, and $\m{W}_o(Y) = \m{O}_Y^\top \m{O}_Y$. The L-ensemble matrix $\m{L}_{\mc{M}} \in \Rn{Np \times Np}$ collects the observability correlations of the measurements
\begin{equation}\label{eq:meas_kernel}
	(\m{L}_{\mc{M}})_{ij} := \tilde{\m{o}}_i^\top\, \tilde{\m{o}}_j, \quad i, j \in \mc{M}.
\end{equation} 

The following theorem considers the $k$-DPP of $\m{L}_{\mc{M}}$, denoted as $\mc{P}_{\m{L}_{\mc{M}}}^{n}$. The DPP samples $n$ measurement rows rather than sensor nodes. A sensor configuration $Y$ is assigned the probability of the event in which every sampled row is generated by the sensors in $Y$. Maximizing this probability over configurations of cardinality $r$ therefore recovers the sensor selection problem under the parameterized Gramian.

\begin{theorem}[\textit{DPP for linear output matrices}]\label{thm:general_C}
	Let $\rank(\m{O}) = n$, so that $\m{W}_o \succ \m{0}$, and let $Y' \sim \mc{P}_{\m{L}_{\mc{M}}}^{n}$ with the kernel~\eqref{eq:meas_kernel}. Then, \textit{(i)} the nonzero eigenvalues of $\m{L}_{\mc{M}}$ are the eigenvalues of $\m{W}_o$, and $e_n(\lambda_1, \ldots, \lambda_{Np}) = \det(\m{W}_o)$; 
	
	\noindent \textit{(ii)} every $S \subseteq \mc{M}$ with $|S| = n$ satisfies
	\begin{equation}\label{eq:proj_prob}
		\mc{P}_{\m{L}_{\mc{M}}}^{n}(S) = \det\big((\m{L}_{\mc{M}})_S\big)/\det(\m{W}_o) = \det(\m{O}_{S})^2/\det(\m{W}_o);
	\end{equation}
	\textit{(iii)} the marginal kernel of $\mc{P}_{\m{L}_{\mc{M}}}^{n}$ is the orthogonal projection
	\begin{equation}\label{eq:proj_kernel}
		\m{K}_o = \m{O}\,\m{W}_o^{-1}\,\m{O}^\top \in \Rn{Np \times Np},
	\end{equation}
	where every $A \subseteq \mc{M}$ satisfies $\prob(A \subseteq Y') = \det((\m{K}_o)_A)$, and the marginal inclusion probability of measurement $i$ is $(\m{K}_o)_{ii} = \tilde{\m{o}}_i^\top \m{W}_o^{-1} \tilde{\m{o}}_i$; and
	
	\noindent \textit{(iv)} every sensor configuration $Y \subseteq \mc{Y}$ satisfies
	\begin{equation}\label{eq:sensor_incl}
		\prob(Y' \subseteq \mc{M}_Y) = \frac{\det(\m{W}_o(Y))}{\det(\m{W}_o)}.
	\end{equation}
\end{theorem}
\begin{proof}
	For \emph{(i)}, the matrices $\m{L}_{\mc{M}} = \m{O}\m{O}^\top$ and $\m{W}_o = \m{O}^\top\m{O}$ share their nonzero eigenvalues~\cite[Sec.~3.3]{Kulesza2012}, and $\rank(\m{L}_{\mc{M}}) = \rank(\m{O}) = n$. In $e_n$, every product of $n$ eigenvalues that contains a zero eigenvalue vanishes, and the remaining product collects the $n$ nonzero eigenvalues, thus $e_n(\lambda_1, \ldots, \lambda_{Np}) = \prod_{i=1}^{n} \lambda_i(\m{W}_o) = \det(\m{W}_o)$. For \emph{(ii)}, \propref{pro:kdpp_map} gives $\mc{P}_{\m{L}_{\mc{M}}}^{n}(S) = \det((\m{L}_{\mc{M}})_S)/e_n$; since $(\m{L}_{\mc{M}})_S = \m{O}_S\m{O}_S^\top$ with $\m{O}_S \in \Rn{n \times n}$, the determinant is $\det(\m{O}_S)^2$, and \emph{(i)} yields~\eqref{eq:proj_prob}. For \emph{(iii)}, the $k$-DPP is a mixture of DPPs over subsets of $k$ eigenvectors, weighted by the products of the corresponding eigenvalues~\cite[Sec.~5.2]{Kulesza2012}. For $k = n = \rank(\m{L}_{\mc{M}})$, the only nonzero weight belongs to the eigenvectors of the $n$ nonzero eigenvalues, so the mixture reduces to the DPP whose marginal kernel is the orthogonal projection onto their span, the column space of $\m{O}$. This projection is~\eqref{eq:proj_kernel}, where $\m{K}_o^\top = \m{K}_o$, $\m{K}_o^2 = \m{O}\m{W}_o^{-1}(\m{O}^\top\m{O})\m{W}_o^{-1}\m{O}^\top = \m{K}_o$, and $\m{K}_o\m{O} = \m{O}$. The sampled set is therefore a DPP with marginal kernel $\m{K}_o$, for which $\prob(A \subseteq Y') = \det((\m{K}_o)_A)$ (\defref{def:marginal_kernel}). Its trace is $\trace(\m{W}_o^{-1}\m{W}_o) = n$, and the diagonal entries of~\eqref{eq:proj_kernel} yield the marginal inclusion probabilities $(\m{K}_o)_{ii} = \tilde{\m{o}}_i^\top \m{W}_o^{-1} \tilde{\m{o}}_i$ with $A = \{i\}$. For \emph{(iv)}, the event $\{Y' \subseteq \mc{M}_Y\}$ is the disjoint union of the events $\{Y' = S\}$ over the subsets $S \subseteq \mc{M}_Y$ with $|S| = n$, so \emph{(ii)} gives
	\begin{equation*}
		\prob(Y' \subseteq \mc{M}_Y)
		= \sum_{\substack{S \subseteq \mc{M}_Y \\ |S| = n}} \frac{\det(\m{O}_{S})^2}{\det(\m{W}_o)}
		= \frac{\det(\m{W}_o(Y))}{\det(\m{W}_o)},
	\end{equation*}
	where the Cauchy-Binet formula reduces the sum of squared minors into $\det(\m{O}_Y^\top \m{O}_Y) = \det(\m{W}_o(Y))$. This concludes the proof.
\end{proof}

\theoref{thm:general_C} extends the relation between DPPs and node selection beyond the setting where a node represents a state variable. The kernel of \theoref{thm:gramian_dpp} is the Gram matrix of the columns of $\m{O}$, and $\m{L}_{\mc{M}}$ is the Gram matrix of its rows; the Gramian $\m{W}_o$ is the dual representation of $\m{L}_{\mc{M}}$~\cite[Sec.~3.3]{Kulesza2012} and carries its nonzero eigenvalues. The ground set therefore indexes the $Np$ measurements, while the normalization and the spectral computations are carried out on the $n \times n$ Gramian. Since $\det(\m{W}_o)$ does not depend on $Y$, maximizing the $\logdett$ Gramian measure $\logdet{\m{W}_o(Y)}$ over sensor configurations of cardinality $r$ also maximizes the probability~\eqref{eq:sensor_incl} that the DPP selects only measurements of the chosen sensors. We note that~\eqref{eq:sensor_incl} is the probability of an event rather than of a single subset; it aggregates the probabilities~\eqref{eq:proj_prob} of all outcomes drawn from the rows of the sensors in $Y$, and the Cauchy-Binet formula reduces this sum of principal minors into $\det(\m{W}_o(Y))$. By the complement property of DPPs (\coref{cor:hodge}), the event probability is itself determinantal, $\prob(Y' \subseteq \mc{M}_Y) = \det(\m{I} - (\m{K}_o)_{\bar{\mc{M}}_Y})$, where $\bar{\mc{M}}_Y := \mc{M} \setminus \mc{M}_Y$ collects the rows of the sensors outside $Y$.

\subsection{Error covariance ratio}\label{app:est_generalC}
The following corollary relates the probability~\eqref{eq:sensor_incl} to the accuracy with which the initial state is estimated from the sensors in $Y$, and provides a lower bound on this probability when the Gramian of the configuration approximates the full Gramian. In contrast to \propref{prop:estimation}, where the $\logdett$ of the principal submatrix $(\m{W}_o)_Y$ is considered, the estimate here considers the parameterized Gramian $\m{W}_o(Y)$. Consider the measurements $\m{y}[k] = \m{C}\m{x}[k] + \m{v}[k]$ over the horizon $N$, where $\m{v}[k]$ is zero mean Gaussian with covariance $\sigma^2 \m{I}$. The least squares estimate of $\m{x}[0]$ from the sensors in $Y$ has error covariance $\m{\Sigma}(Y) := \sigma^2 \m{W}_o(Y)^{-1}$, and $\m{\Sigma}(\mc{Y}) = \sigma^2 \m{W}_o^{-1}$ when all $p$ sensors are measured.

\begin{mycor}[\textit{Error covariance ratio for linear output matrices}]\label{cor:generalC_est}
	Let $Y' \sim \mc{P}_{\m{L}_{\mc{M}}}^{n}$, and let $\m{W}_o(Y) \succ \m{0}$ for a sensor configuration $Y \subseteq \mc{Y}$. Then
	\begin{equation}\label{eq:est_ratio}
		\prob(Y' \subseteq \mc{M}_Y) = \frac{\det(\m{\Sigma}(\mc{Y}))}{\det(\m{\Sigma}(Y))},
	\end{equation}
	and if $(1-\epsilon)\, \m{W}_o \preceq \m{W}_o(Y)$ for some $\epsilon \in [0,1)$, then
	\begin{equation}\label{eq:est_bound}
		\prob(Y' \subseteq \mc{M}_Y) \geq (1-\epsilon)^n, \;\;\; \det(\m{\Sigma}(Y)) \leq \frac{\det(\m{\Sigma}(\mc{Y}))}{(1-\epsilon)^{n}}.
	\end{equation}
\end{mycor}
\begin{proof}
	The proof follows from stacking the measurements of the sensors in $Y$ over the horizon $N$, which yields $\m{Y} = \m{O}_Y \m{x}[0] + \bm{\eta}$; the least squares estimate of $\m{x}[0]$ has error covariance $\m{\Sigma}(Y) = \sigma^2 (\m{O}_Y^\top \m{O}_Y)^{-1} = \sigma^2 \m{W}_o(Y)^{-1}$, and the case $Y = \mc{Y}$ gives $\m{\Sigma}(\mc{Y}) = \sigma^2 \m{W}_o^{-1}$. It follows that $\det(\m{\Sigma}(\mc{Y}))/\det(\m{\Sigma}(Y)) = \det(\m{W}_o(Y))/\det(\m{W}_o)$, which equals $\prob(Y' \subseteq \mc{M}_Y)$ by~\eqref{eq:sensor_incl}; this yields~\eqref{eq:est_ratio}. For~\eqref{eq:est_bound}, the ordering $(1-\epsilon)\m{W}_o \preceq \m{W}_o(Y)$ implies $(1-\epsilon)^n \det(\m{W}_o) \leq \det(\m{W}_o(Y))$, and thus $\prob(Y' \subseteq \mc{M}_Y) \geq (1-\epsilon)^n$ by~\eqref{eq:sensor_incl}; the covariance bound then follows from~\eqref{eq:est_ratio}. This concludes the proof.
\end{proof}

\coref{cor:generalC_est} shows that the probability~\eqref{eq:sensor_incl} equals the ratio between the determinant of the estimation error covariance when all sensors are measured. Given that the determinant of the error covariance quantifies the volume of the state estimation confidence ellipsoid~\cite{Joshi2009}, sensor configurations with probability close to one estimate $\m{x}[0]$ nearly as accurately as the full sensor set. The condition $(1-\epsilon)\m{W}_o \preceq \m{W}_o(Y)$ in~\eqref{eq:est_bound} is the one sided form of the Gramian approximation that defines the sparse actuator schedules of~\cite{Siami2021}. This bound~\eqref{eq:est_bound} translates this Gramian approximation guarantee into a lower bound on the DPP probability of the configuration and on its estimation accuracy. Since the exponent in~\eqref{eq:est_bound} grows with the state dimension; it is informative in its normalized form $\prob(Y' \subseteq \mc{M}_Y)^{1/n} \geq 1-\epsilon$. This is consistent with the normalized determinant measure $\det(\cdot)^{1/n}$ used to compare Gramians in~\cite{Siami2021}.

\balance
\bibliographystyle{IEEEtran}
\bibliography{library_cdc2026}

\end{document}